\documentclass[11pt, letterpaper]{article}
\usepackage{xcolor}
\usepackage{colortbl}
\usepackage{verbatim}
\usepackage{enumerate}
\usepackage{natbib}

\usepackage{fancyhdr}
\usepackage{amsmath, amsfonts}
\usepackage{natbib}
\usepackage{amssymb}
\usepackage{latexsym}
\usepackage{theorem}
\usepackage{enumerate}
\usepackage{graphicx}
\usepackage{rotating}
\usepackage{pdfpages}

\newtheorem{defn}{Definition}[section]
\newtheorem{lemma}[defn]{Lemma}
\newtheorem{theorem}[defn]{Theorem}

\newtheorem{assumption}[defn]{Assumption}

\newtheorem{estprocedure}[defn]{Estimation Procedure}

\newcommand{\cip}{\mbox{$\perp\!\!\!\perp$}}

\title{How estimating nuisance parameters can reduce the variance (with consistent variance estimation)}
\author{Judith J. Lok\\
Department of Mathematics and Statistics\\
Boston University\\
jjlok@bu.edu}

\begin{document}

\maketitle

\section*{Abstract}

We often estimate a parameter of interest $\psi$ when the identifying conditions involve a nuisance parameter $\theta$. Examples from causal inference are Inverse Probability Weighting, Marginal Structural Models and Structural Nested Models. To estimate treatment effects from observational data, these methods posit a (pooled) logistic regression model for the treatment (initiation) and/or censoring probabilities, and estimate these in a first step. These methods can all be shown to be based on unbiased estimating equations. First, we provide a general formula for the variance of the parameter of interest $\psi$ when the nuisance parameter $\theta$ is estimated in a first step, using the Partition Inverse Formula. Then, we present 4 results for estimators $\hat{\psi}$ based on unbiased estimating equations including a nuisance parameter $\theta$ which is estimated by solving (partial) score equations, in settings where $\psi$ does not depend on $\theta$. This regularly happens in causal inference if $\theta$ describes the treatment allocation probabilities, in missing data settings where $\theta$ describes the missingness probabilities, and in measurement error settings  where $\theta$ describes the measurement error distribution. 1. Counter-intuitively, the limiting variance of $\hat{\psi}$ is typically smaller when $\theta$ is estimated, compared to if a known $\theta$ were plugged in. 2. If estimating $\theta$ is ignored, the resulting sandwich estimator for the variance of $\hat{\psi}$ is conservative. 3. A consistent estimator for the variance of $\hat{\psi}$ provides results fast: no bootstrap. 4. If $\hat{\psi}$ with the true $\theta$ plugged in is efficient, the limiting variance of $\hat{\psi}$ does not depend on whether $\theta$ is estimated, and the sandwich estimator for the variance of $\hat{\psi}$ ignoring estimation of $\theta$ is consistent. To illustrate we use observational data to estimate 1. the effect of caz-avi versus colistin on clinical outcomes in patients with resistant bacterial infections and 2. how the effect of one year of antiretroviral treatment depends on its initiation time in HIV-infected patients.

\section{Introduction}

In many applications, to estimate a parameter $\psi$ of interest, a finite-dimensional nuisance parameter $\theta$ is estimated in a first step. Examples include Inverse Probability Weighting (\cite{RRZ}, \cite{Bang}), Marginal Structural Models (\cite{MSM1}, \cite{MSM2}), and Structural Nested Models (\cite{Enc}, \cite{SNart}). Confidence interval/variance estimation for the parameter of interest $\psi$ in such settings is complicated, since the variance of the resulting estimator $\hat{\psi}$ may depend on that the nuisance parameter $\theta$ is estimated, and how. This often leads analysts to use the bootstrap (see for example \cite{MSM1} and \cite{MSM2}). The bootstrap (\cite{Efron}) is often successfully used for a variety of estimation problems, but especially for larger datasets and complicated analyses the bootstrap can require considerable computing time.

Existing theory for consistent estimation of the variance in settings with nuisance parameters includes but is not restricted to the following. For longitudinal studies analyzed with Inverse Probability of Censoring Weighting, \cite{RR95} provided variance estimates. For Inverse Probability of Treatment Weighting for point treatment, \cite{lunceford2004stratification} provided variance estimates. \cite{lunceford2004stratification} notice that estimating the nuisance parameters does not affect the variance for the doubly robust alternative of \cite{RRZ94} (also described in \citeauthor{Bang} \citeyear{Bang}). For continuous-time Structural Nested Models, \cite{Lok} and \cite{ASarx} provided variance estimates. \cite{henmi2004paradox} focused on settings where the parameter of interest $\psi$ and the nuisance parameter $\theta$ appear in a joint likelihood and vary independently, and where the scores with respect to $\psi$ and the parametric part of the nuisance parameter $\theta$ are both well-defined. Their results in this setting follow projections and geometric arguments with respect to the scores with respect to both $\psi$ and $\theta$. Their setting differs from the examples in this article, where the parameter of interest $\psi$ does not parameterize the likelihood; instead, it is a function of the finite-dimensional and infinite-dimensional parameters involved in the likelihood, even in the well-studied setting of estimating an overall mean based on observational data using Inverse Probability of Treatment Weighting (Sections~\ref{IPTWpoint} and~\ref{cazavi}; see also Section~\ref{Discussion}).

Various researchers have worked on variance estimation for a propensity score-weighted marginal hazard ratio depending on treated/untreated at baseline. \cite{Hajage} used an influence function linearization technique. \cite{shu2021variance} developed an asymptotically equivalent closed form correction to the sandwich estimator in this setting, and proved that without this correction, the sandwich estimator is conservative. \cite{Austin} carried out simulation studies in this setting and showed that the sandwich estimator, also described in \cite{Joffe2004model}, can also be conservative in practice. 

Inverse Probability Weighting, Marginal Structural Models, and Structural Nested Models can all be shown to be based on estimating equations that would be unbiased were the nuisance parameter $\theta$ known. The current paper includes all settings with independent, identically distributed observations where the parameter of interest $\psi$ is estimated based on estimating equations that would be unbiased were the nuisance parameter $\theta$ known, and where $\theta$ is estimated by solving unbiased estimating equations. We consider as a special and fairly common case settings where $\theta$ is estimated by solving (partial) score equations and $\psi$ does not depend on $\theta$. Causal inference includes many such settings: often, the propensity scores (that is, the probabilities of treatment given baseline covariates, or the probabilities of a time-dependent treatment given the past) are modeled with (pooled) logistic regression depending on a nuisance parameter $\theta$. $\theta$ is then estimated with standard software, which solves (partial) score equations.

For longitudinal studies analyzed with Inverse Probability of Censoring Weighting, our variance derivations lead to the same results as in \cite{RR95}. For Inverse Probability of Treatment Weighting for point treatment, our variance derivations lead to the same results as in \cite{lunceford2004stratification}. We thus generalize these results to many other estimation settings with nuisance parameters.

Section~\ref{cazavi} illustrates the approach by estimating the effect of caz-avi versus colistin for the treatment of resistant bacterial infections. The data are observational, and the analysis involves Inverse Probability of Treatment Weighting for point treatment; the original paper \cite{Duin} used the bootstrap to create confidence intervals. Section~\ref{simulations} describes simulations modeled after Section~\ref{cazavi}. Section~\ref{ART} illustrates the approach by estimating how the effect of one year of ART depends on its initiation time in HIV-infected patients. The data are observational, and the analysis involves coarse Structural Nested Mean Models; the original paper \cite{Victor} used the bootstrap to create confidence intervals, which took over 8 hours on a large computing cluster; the sandwich estimator took a few minutes on a laptop.

We expect our approach to be highly relevant for applied research inside and outside of causal inference, where one often needs to estimate a nuisance parameter $\theta$ before being able to estimate the parameter of interest $\psi$.

\section{Setting and notation}\label{setting}

Throughout this paper we assume independent, identically distributed observations $(Y_i,Z_i)$, $i=1,\ldots,n$. Interest lies in a finite-dimensional parameter $\psi$, and we have some set of unbiased estimating equations for the true parameter of interest $\psi$ were a nuisance parameter $\theta$ known,
\begin{equation*}
\frac{1}{n}\sum_{i=1}^n U_{1}(Y_i,Z_i;\psi,\theta)=0,
\end{equation*}
or equivalently, with $P_n$ the empirical average over $i=1,\ldots, n$,
\begin{equation*}
P_n U_{1}(\psi,\theta)=0.
\end{equation*}
Throughout, we assume
\begin{assumption}\label{uEEs} {\bf Unbiased estimating equations.}
\begin{equation*}
EU_{1}(Y_i,Z_i;\psi^*,\theta^*)=0
\end{equation*}
for the true parameters $\psi^*$ and $\theta^*$.
\end{assumption}
In many settings, Assumption~\ref{uEEs} follows from model assumptions and identifiability assumptions, as in \cite{MSM1} and \cite{MSM2} 
for Marginal Structural Models, in \cite{Aids} and \cite{ASarx} for Structural Nested Models, and in both our illustrative examples from Section~\ref{Examples}: Inverse Probability of Treatment Weighting for point treatment and time-dependent coarse Structural Nested Mean Models. After finding such identifying result, we will henceforth focus on estimating the zero $\psi^*$ of $EU_1(\psi,\theta^*)=0$, without making assumptions other than the parameterization that identifies $\theta$. Thus, we simply define 

\begin{defn} $\psi^*$ is the solution to $EU_1(\psi,\theta^*)=0$.
\end{defn}
Since the true $\theta$ is unknown, we cannot immediately use $U_{1}$ to estimate $\psi$.

To use $U_1$ to estimate $\psi$, $\theta$ is estimated first. Then, the estimate $\hat{\theta}$ is plugged into the estimating equations for $\psi$ to solve $P_n U_1(\psi,\hat{\theta})=0$:
\begin{estprocedure}\label{proc} $\psi$ is estimated by estimating the nuisance parameter $\theta$ first, and $\hat{\psi}$ is the solution to 
\begin{equation}\label{EE}
\frac{1}{n}\sum_{i=1}^n U_{1}(Y_i,Z_i;\psi,\hat{\theta})=0.
\end{equation}
\end{estprocedure}
Throughout, we assume that $\theta$ is estimated by solving unbiased estimating equations:
\begin{assumption}\label{U2} {\bf Unbiased estimating equations for $\theta$.}
$\theta$ is estimated by solving 
\begin{equation*}
\frac{1}{n}\sum_{i=1}^n U_{2}(Z_i;\theta)=0, \;\;\;\;\;\text{that is,}\;\;\;\;\; P_nU_2(\theta)=0,
\end{equation*}
with
\begin{equation*}
EU_{2}(Z_i;\theta^*)=0.
\end{equation*}
\end{assumption}

Often, we assume a stronger version of Assumption~\ref{U2}, as follows. In general, when we consider the data generating mechanism, we assume that the data have a distribution that is partly parameterized by $\theta$. The concept of partial likelihood was introduced in \cite{pCox}. One orders the parts of the observation $(Z_i,Y_i)$ from a unit $i$ in some way and parameterizes the distribution of some (but not all) of these observation parts given all the ``past'' observation parts.
For example, one could parameterize the distribution $f_\theta(Z_i)$ of the $Z_i$, but not the conditional distribution of the $Y_i$ given $Z_i$. For ease of understanding, in the following one can think of $\theta$ parameterizing the distribution of the $Z_i$ and of the conditional distribution of $Y_i$ given $Z_i$ remaining unspecified, although many applications will be more like the examples in Section~\ref{Examples}, where only part of the distribution of the $Z_i$ is determined by $\theta$. In a large part of our exposition, we will assume the following stronger version of Assumption~\ref{U2}:

\begin{assumption}\label{MPLE}{\bf (Partial) score equations for $\theta$.}
$\theta$ is estimated by maximum (partial) likelihood on (part of) the $Z_i$, and this leads to $\hat{\theta}$ solving (partial) score equations
\begin{equation*}
\frac{1}{n}\sum_{i=1}^n U_{2}(Z_i;\theta)=\frac{1}{n}\sum_{i=1}^n\frac{\partial}{\partial \theta} \log f_\theta(Z_i)=0
\end{equation*}
with
\begin{equation*}
EU_{2}(Z_i;\theta^*)=E\left(\left.\frac{\partial}{\partial \theta}\right|_{\theta^*} \log f_\theta(Z_i)\right)=0.
\end{equation*}
\end{assumption}

As in \cite{Victor}, it is easy to see that Estimation Procedure~\ref{proc}
is equivalent to solving
\begin{equation}\label{stack} P_n\left(\begin{array}{c} U_1(\psi,\theta)\\U_2(\theta)\end{array}\right)=0
\end{equation}
for $(\psi,\theta)$ to obtain $(\hat{\psi},\hat{\theta})$, since the second set of estimating equations $P_n U_2(\theta)=0$ determine $\hat{\theta}$. One would not necessarily want to actually calculate $(\hat{\psi},\hat{\theta})$ using equation~(\ref{stack}), since standard software can often be used to obtain $\hat{\theta}$ and then other standard software can often be used to obtain $\hat{\psi}$ after plugging in $\hat{\theta}$. However, the stacked equations~(\ref{stack}) are very useful for theory development. For example, consistency and asymptotic normality of $\hat{\psi}$ follows, under regularity conditions, from \cite{Vaart} Theorems~5.9 and~5.21, since (\ref{stack}) forms unbiased estimating equations for $(\psi,\theta)$ jointly. 

For the result that the sandwich estimator ignoring estimation of $\theta$ leads to conservative inference, we also need:
\begin{assumption}\label{psinotftheta} The parameter of interest $\psi$ does not depend on the nuisance parameter $\theta$.
\end{assumption}
To clarify what is meant by the parameter $\psi$ not depending on $\theta$, notice that Assumption~\ref{psinotftheta} is used together with Assumption~\ref{MPLE}. Under Assumption~\ref{MPLE}, the likelihood $f(y,z)$ depends on $\theta$ and a (possibly infinite dimensional) parameter say $\eta$, and $f(y,z)$ factorizes as $f(y,z)=f_{\theta}(z)f_\eta(y,z)$. Assumption~\ref{psinotftheta} then states that $\psi$, a function of the distribution of $(y,z)$, is a function of $\eta$ and not $\theta$. Settings where Assumption~\ref{psinotftheta} often holds are, for example, 1.\ observational studies where $\theta$ describes the treatment allocation, if interest lies in overall treatment effects or conditional treatment effects, 2.\ settings with missing data where $\theta$ describes the missingness probabilities, and 3.\ settings with measurement error where $\theta$ describes the measurement error distribution. Assumption~\ref{psinotftheta} does not hold for example for the treatment effect in the treated.

\section{Examples}\label{Examples}

We provide some illustrative causal inference examples.

\subsection{Example: estimating treatment effects from observational data}\label{IPTWpoint}

Consider a setting with $i$ indexing patient $i$, $L_i$ baseline covariates, $A_i$ a binary treatment, and $Y_i$ the outcome of interest of patient $i$. Suppose treatment $A_i$ is not randomized. Here, $(Y_i,Z_i)=(Y_i,A_i,L_i)$. 
Write $Y_i^{(1)}$ for the outcome of patient $i$ under treatment $a=1$, and $Y_i^{(0)}$ for the outcome of patient $i$ under ``no treatment'' $a=0$. These are so-called counterfactual outcomes (\citeauthor{prop}, \citeyear{prop}), since we can only observe one of these for each patient, leaving either $Y_i^{(1)}$ or $Y_i^{(0)}$ counter-to-fact. Of interest are often the expected outcome under treatment $\psi_1^*=E\bigl(Y^{(1)}\bigr)$, the expected outcome under no treatment $\psi_0^*=E\bigl(Y^{(0)}\bigr)$, and their difference
\begin{equation}\label{trteffect}
E\bigl(Y_i^{(1)}-Y_i^{(0)}\bigr)=\psi_1^*-\psi_0^*.
\end{equation}

Suppose there is ``No Unmeasured Confounding'': with $\cip$ meaning ``is independent of'' (\cite{Dawid}):
 \begin{equation}\label{nuc}
 Y_i^{(1)}\cip A_i \mid L_i \;\;\;\text{  and  }\;\;\; Y_i^{(0)}\cip A_i \mid L_i
\end{equation}
(see, for example, \citeauthor{prop} \citeyear{prop}). 
The interpretation of this Assumption of No Unmeasured Confounding is that doctors and patients may have used the observed baseline covariates $L_i$ for their treatment decisions $A_i$, but other than that, a patient's prognosis (represented by $Y_i^{(1)}$ and $Y_i^{(0)}$) did not affect their treatment decisions $A_i$. That is, treatment decisions did not peek into the future, $Y_i^{(1)}$ or $Y_i^{(0)}$, based on information other than observed in $L_i$. For example, if patients' blood pressure or CD4 count affected the decisions to treat or not, blood pressure or the CD4 count needs to be collected in the $L_i$.

Under No Unmeasured Confounding, causal inference provides various ways to estimate the average treatment effect $\psi_1-\psi_0$. Inverse Probability of Treatment Weighting (IPTW, see e.g. \citeauthor{MSM1} \citeyear{MSM2}, \citeauthor{MSM2} \citeyear{MSM2}, and \citeauthor{Bang} \citeyear{Bang}) first estimates how the treatment $A_i$ depends on the baseline covariates $L_i$, typically by specifying a parametric model $p_\theta(A_i=1|L_i)$ for $P(A_i=1|L_i)$, the so-called propensity score (\citeauthor{prop}, \citeyear{prop}). This leads to 
\begin{equation*}
f_{\theta,\cdot}(Y,A,L)=f(L)p_\theta(A|L) f(Y|A,L),
\end{equation*}
with the distribution of $L$ and the distribution of $Y$ given $(A,L)$ remaining completely unspecified. 

To estimate $\theta$, maximum partial likelihood maximizes the likelihood
\begin{equation*}
\prod_{i=1}^n f(L_i)p_\theta(A_i|L_i) f(Y_i|A_i,L_i)
\end{equation*}
over $\theta$, which is the same as maximizing the partial likelihood
\begin{equation*}
\prod_{i=1}^n p_\theta(A_i|L_i).
\end{equation*}
For most models $p_\theta$ this leads to solving partial score equations
\begin{equation*}
\frac{1}{n}\sum_{i=1}^n \frac{\partial}{\partial \theta} \log p_\theta(A_i|L_i)=0.
\end{equation*}
Here, $U_{2}(Z_i;\theta)=\frac{\partial}{\partial \theta} \log p_\theta(A_i|L_i)$.
$p_\theta(A_i=1|L_i)$ is often modeled with a logistic regression model. For a logistic regression model $p_\theta(A_i=1|L_i)=e^{\theta L_i}/(1+e^{\theta L_i})$, with $\theta$ and $L_i$ say $m$-dimensional and $L_i$ usually including $1$ in the first row to allow for an intercept, this leads to the familiar score equations for logistic regression,
\begin{equation*}
\frac{1}{n}\sum_{i=1}^n  L_i (A_i-p_\theta(A_i=1|L_i))=0.
\end{equation*}
In this case, the partial scores for $\theta$ are
\begin{equation}\label{IPTWpointU2}
 U_{2}(Z_i;\theta)=L_i (A_i-p_\theta(A_i=1|L_i)).   
\end{equation}

IPTW then proceeds by estimating the mean $\psi_1$ under treatment ($a=1$) and the mean $\psi_0$ under no treatment ($a=0$) by solving
\begin{equation*}
\frac{1}{n}\sum_{i=1}^n \frac{1_{A_i=1}}{P_{\hat{\theta}}\left(A_i=1|L_i\right)}\left(Y_i-\psi_1\right)=0,\hspace*{0.5cm}
\frac{1}{n}\sum_{i=1}^n \frac{1_{A_i=0}}{P_{\hat{\theta}}\left(A_i=0|L_i\right)}\left(Y_i-\psi_0\right)=0
\end{equation*}
to obtain $\hat{\psi}_1$ and $\hat{\psi}_0$ and their difference $\hat{\psi}_1-\hat{\psi}_0$, the estimated effect of treatment $A$ on the outcome $Y$. To avoid division by $0$ one needs positivity (\citeauthor{prop} \citeyear{prop}): $P_{\theta}\left(A=1|L\right)$ needs to be bounded away from $0$ and $1$. The intuition behind IPTW is that when we want to estimate say $EY_i^{(1)}$ and e.g.\ only $1/3$ of patients with a specific value of $L$ receive treatment $a=1$, those patients need to represent 3 patients each: themselves and 2 others who were not taking treatment $a=1$. Hence, they need to carry a weight of one over their probability to receive $a=1$.

In this example, $U_1$ can be seen to lead to unbiased estimating equations as follows, for $\psi_1$ (and similarly for $\psi_0$):
\begin{eqnarray*}
E\left(U_1(\psi_1^*,\theta^*)\right)
&=&E\left(\frac{1_{A_i=1}}{P_{\theta^*}\left(A_i=1|L_i\right)}\left(Y_i^{(1)}-\psi_1^*\right)\right)\\
&=&E\left(E\left[\frac{1_{A_i=1}}{P_{\theta^*}\left(A_i=1|L_i\right)}\left(Y_i^{(1)}-\psi_1^*\right)\mid Y_i^{(1)},L_i\right]\right)\\
&=&E\left(\frac{E\left[1_{A_i=1}\mid Y_i^{(1)},L_i\right]}{P_{\theta^*}\left(A_i=1|L_i\right)}\left(Y_i^{(1)}-\psi_1^*\right)\right)\\
&=&E\left(\frac{E\left[1_{A_i=1}\mid L_i\right]}{P_{\theta^*}\left(A_i=1|L_i\right)}\left(Y_i^{(1)}-\psi_1^*\right)\right)\\
&=&E\left(Y_i^{(1)}-\psi_1^*\right)=0,
\end{eqnarray*}
where for the first equality we use that on $A=1$, we observe the outcome $Y^{(1)}$ under treatment, for the second equality we use the Law of Iterated Expectations, and for the fourth equality we use the Assumption of No Unmeasured Confounding (\ref{nuc}).

\subsection{Example: estimating treatment effects from longitudinal observational data} \label{long}


Consider an observational longitudinal study with $i$ indexing patient $i$, $L_{ki}$ the covariates at time $k$, $A_{ki}$ the binary treatment starting at time $k$ and lasting until just before time $k+1$, and $Y_i$ the outcome at the end of the study. The observed data are
\begin{equation*}
(L_{0i},A_{0i},L_{1i},A_{1i},\ldots,L_{Ki},A_{Ki},Y_i)    
\end{equation*}
for $i=1\ldots n$ patients. Here $Z_i=(L_{0i},A_{0i},L_{1i},A_{1i},\ldots,L_{Ki},A_{Ki})$, and $\theta$ could parameterize the conditional distribution of the $A_{ki}$ given the past,
\begin{equation*}p_\theta(A_{ki}|L_{0i},A_{0i},L_{1i},A_{1i},\ldots,A_{k-1i},L_{ki}),
\end{equation*}
leaving the conditional distributions of the $L_{ki}$ and $Y_i$ unspecified. With $\overline{L}_k=(L_0,L_1,\ldots,L_k)$ and $\overline{A}_k=(A_0,A_1,\ldots,A_k)$ the covariate- and treatment histories, this leads to 
\begin{equation}\label{density}
f_{\theta,\cdot}(Y,Z)=\prod_{k=0}^K p_\theta(A_k|\overline{L}_k,\overline{A}_{k-1}) f(L_k|\overline{L}_{k-1},\overline{A}_{k-1})f(Y|\overline{L}_{K},\overline{A}_{K}),
\end{equation}
with the $f(L_k|\overline{L}_{k-1},\overline{A}_{k-1})$ and $f(Y|\overline{L}_{K},\overline{A}_{K})$ unspecified.
Thus, we just model the probability/likelihood of treatment, ignoring how the previous information came about (hence ``partial'' likelihood). Equation~(\ref{density}) is the common starting point for Marginal Structural Models (\cite{MSM1}, \cite{MSM2}) and Structural Nested Models (\cite{Aids}, \cite{smoke}, \cite{SNart}).

To estimate $\theta$, maximum partial likelihood maximizes the likelihood
\begin{equation*}
\prod_{i=1}^n\prod_{k=0}^K p_\theta(A_{ki}|\overline{L}_{ki},\overline{A}_{(k-1)i}) f(L_{ki}|\overline{L}_{(k-1)i},\overline{A}_{(k-1)i})f(Y_i|\overline{L}_{Ki},\overline{A}_{Ki})
\end{equation*}
over $\theta$, which is the same as maximizing the partial likelihood over $\theta$:
\begin{equation*}
\prod_{i=1}^n \prod_{k=0}^K p_\theta(A_{ki}|\overline{L}_{ki},\overline{A}_{(k-1)i}).
\end{equation*}
For most models $p_\theta$ this leads to solving partial score equations
\begin{equation*}
\frac{1}{n}\sum_{i=1}^n \sum_{k=0}^K \frac{\partial}{\partial \theta}p_\theta(A_{ki}|\overline{L}_{ki},\overline{A}_{(k-1)i})=0.
\end{equation*}
For a pooled logistic regression model such as $p_\theta(A_{ki}=1|\overline{L}_{ki},\overline{A}_{k-1i})=e^{\theta_1 L_{ki}+\theta_2 A_{k-1i}}/(1+e^{\theta_1 L_{ki}+\theta_2 A_{k-1i}})$, with $\theta_1$ and $L_{ki}$ say $m$-dimensional and $L_{ki}$ usually including $1$ in the first row to allow for an intercept, this leads to the familiar score equations for logistic regression,
\begin{equation}\label{scoreMSM}
\frac{1}{n}\sum_{i=1}^n \sum_{k=0}^K \left(\begin{array}{c} L_{ki}\\ A_{(k-1)i} \end{array}\right) (A_{ki}-p_\theta(A_{ki}=1|\overline{L}_{ki},\overline{A}_{(k-1)i}))=0.
\end{equation}
In this case, 
\begin{equation*}
U_{2}(Z_i;\theta)=\sum_{k=0}^K \left(\begin{array}{c} L_{ki}\\ A_{k-1i} \end{array}\right) (A_{ki}-p_\theta(A_{ki}=1|\overline{L}_{ki},\overline{A}_{k-1i}))
\end{equation*}
are the partial scores for $\theta$.

Many causal inference methods to estimate treatment effects from observational data assume the following version of No Unmeasured Confounding:
\begin{assumption}\label{nucwithtime} {\bf No Unmeasured Confounding.}
\begin{equation*}A_{ki}\cip (\overline{L}_{Ki}^{(\overline{a}_K)},Y_i^{(\overline{a}_K)}) \mid \overline{L}_{ki},\overline{A}_{k-1i}=\overline{a}_{k-1},
\end{equation*}
where $\overline{L}_{Ki}^{(\overline{a}_K)}$ and $Y_i^{(\overline{a}_K)}$ are the (counterfactual) covariate history and outcome of patient $i$ under treatment regime $\overline{a}_K$.
\end{assumption}
The interpretation of this Assumption of No Unmeasured Confounding is the same as for point treatment, see equation~(\ref{nuc}): we assume that treatment decisions may depend on a patient's past observed information, and not further on a patient's prognosis.

\subsubsection{Example: Inverse Probability of Treatment Weighting for time-dependent treatments}\label{IPTWtimedep}

For time-dependent treatments, Inverse Probability of Treatment Weighting (IPTW) often posits a parametric model for the mean outcome under any fixed treatment regimen $\overline{a}_K$, say
\begin{equation*}
E\left(Y_i^{(\overline{a}_K)}\right)=\mu_\psi^{(\overline{a}_K)},
\end{equation*}
or, conditional on baseline covariates $X$,
\begin{equation*}
E\left[Y_i^{(\overline{a}_K)}\mid X_i\right] =\mu_\psi^{(\overline{a}_K)}(X_i);
\end{equation*}
for example, with $\psi_0$ typically a vector,
\begin{equation*}
E\left[Y_i^{(\overline{a}_K)}\mid X_i\right] =\psi_0X_i+\psi_1\sum_{k=0}^Ka_k.
\end{equation*}
IPTW estimates $\psi$ by using
\begin{equation*}
U_1(Y_i,\overline{L}_{Ki},\overline{A}_{Ki};\psi,\theta)=\sum_{\overline{a}_K} \frac{1_{\overline{A}_{Ki}=\overline{a}_{Ki}}}{\prod_{k=0}^K P_{\theta}\left(A_{ki}=a_{ki}|\overline{L}_{ki},\overline{A}_{k-1i}=\overline{a}_{k-1}\right)}\left(Y_i-\mu_\psi^{(\overline{a}_K)}\right),
\end{equation*}
estimating $\theta$ in a first step as in Estimation Procedure~\ref{proc}.
Alternatively, to increase precision, stabilized weights can be used. For the model $\mu_\psi^{(\overline{a}_K)}(X_i)$, stabilized weights add $\prod_{k=0}^K P_{\hat{\theta}}\left(A_{ki}=a_{ki}|X_i,\overline{A}_{k-1i}=\overline{a}_{k-1}\right)$ to the numerator; for $\mu_\psi^{(\overline{a}_K)}$ not conditional on $X_i$, stabilized weights add $\prod_{k=0}^K P_{\hat{\theta}}\left(A_{ki}=a_{ki}|\overline{A}_{k-1i}=\overline{a}_{k-1}\right)$ to the numerator. Similar to the point-treatment setting, it can be shown (see Web-appendix~\ref{Proofs}) that under No Unmeasured Confounding and Positivity,
\begin{equation}\label{IPTWwithtime}
EU_1(Y_i,\overline{L}_{Ki},\overline{A}_{Ki};\psi^*,\theta^*)=0
\end{equation}
for these methods. 

\subsubsection{Example: time-dependent coarse Structural Nested Mean Models (coarse SNMMs)}\label{coarse}

Our explanations for time-dependent coarse Structural Nested Mean Models (coarse SNMMs) are relatively brief; more background and details are provided in \cite{Victor}. Also coarse SNMMs estimate how treatment affects the outcome of interest, from non-randomized data. 

The notation in Section~\ref{coarse} is the same as in Section~\ref{long}. Here, time-dependent outcomes of interest $Y_k$ are measured at times $0,\ldots,K+1$. Treatment can be initiated at any time point, and we are interested in the effect of treatment initiation regardless of whether treatment is continued later (compare with intention-to-treat). $T$ is the time of treatment initiation, with $T_i=K+1$ if patient $i$ does not initiate treatment before time $K+1$.

Write $Y^{(\emptyset)}_{k}$ for the (possibly counterfactual, or not-observed) outcome at
time $k$ without treatment, and $Y^{(m)}_{k}$ for the (possibly counterfactual) outcome at
time $k$ with treatment initiated at time $m$. Time-dependent coarse Structural Nested Models (SNMMs) describe contrasts between outcomes $Y^{(\emptyset)}_{k}$  and $Y^{(m)}_{k}$:

\begin{defn} \label{gamma} {\bf Time-dependent coarse SNMM.} For $k=m,\ldots, K+1$,
\begin{equation*}
\hspace*{3cm}\gamma_k^{(m)}\left(\overline{l}_m\right)=E\left[Y_{k}^{(m)}-Y_{k}^{(\emptyset)}\mid \overline{L}^{(\emptyset)}_{m}=\overline{l}_m,T=m\right].
\end{equation*}
\end{defn}
$\gamma_k^{(m)}\left(\overline{l}_m\right)$ is the expected difference of 
the outcome at time $k$ had the patient initiated treatment at time
$m$, and the outcome at time $k$ had the
patient never initiated treatment. It is the effect of treatment between times $m$ and $k$,  given a patient's pre-treatment information $\overline{l}_m$.

\begin{assumption}{\bf Model.} We assume a correctly specified parametric model $\gamma_\psi$ for $\gamma$.
\end{assumption}
The recurring example in \cite{Victor} has
\begin{equation}\label{modelex}
\gamma^{(m)}_{k,\psi}\left(\overline{l}_m\right)=
\left(\psi_1+\psi_2m+\psi_3m^2\right)(k-m)1_{\left\{k>m\right\}},
\end{equation} with $(k-m)$ the treatment duration from
month $m$ to month $k$.

As in \cite{Aids}, \cite{Enc}, \cite{MSM1}, \cite{SNart},
\cite{Rnoncomp}, \cite{Victor}, and \cite{optLok}, to distinguish between the reason to treat and the treatment's effect, coarse SNMMs need all patient-level information that both: (1) influences treatment decisions
and (2) possibly predicts a patient's prognosis with respect to the
outcome of interest. $Y^{(\emptyset)}_{k}$, the outcome at
time $k$ without treatment, reflects a patient's
prognosis with respect to the outcome of interest. No Unmeasured
Confounding means that treatment decisions at time $m$ ($A_m$) should be
independent of this (not always observed) prognosis $Y^{(\emptyset)}_{k}$ ($k>m$), given the observed
past treatment and covariate history $(\overline{L}_{m},\overline{A}_{m-1})$:\medskip

\begin{assumption} \label{intconf}{\bf (No Unmeasured Confounding -
    formalization)}. $A_{mi}\cip Y^{(\emptyset)}_{ki}
\mid \overline{L}_{mi},\overline{A}_{m-1i}=\overline{0}$
for $k>m$.
\end{assumption}
For the parameter of interest $\psi$ to satisfy Assumption~\ref{psinotftheta}, we adopt the following stronger version of No Unmeasured Confounding, albeit with the same interpretation as Assumption~\ref{intconf}:\medskip

\begin{assumption} \label{intconf2}{\bf (No Unmeasured Confounding -
    formalization)}. $A_{pi}\cip (Y^{(\emptyset)}_{ki},Y_{ki}^{(m)},\overline{L}_{mi}^{(\emptyset)})
\mid \overline{L}_{pi},\overline{A}_{p-1i}=\overline{0}$
for $k>m$ and $p\leq m$.
\end{assumption}
Under Assumption~\ref{intconf2}, it can be shown that
\begin{equation*}
\hspace*{3cm}\gamma_k^{(m)}\left(\overline{l}_m\right)=E\left[Y_{k}^{(m)}-Y_{k}^{(\emptyset)}\mid \overline{L}^{(\emptyset)}_{m}=\overline{l}_m\right],
\end{equation*}
so that Assumption~\ref{psinotftheta} is satisfied.

If a patient is not treated until time $k$, there is no difference in
treatment between $Y_{k}$ and $Y^{(\emptyset)}_{k}$, so then it is reasonable to assume that until time $k$, the observed outcomes and covariates equal those without treatment:

\begin{assumption} \label{cons}\emph{(Consistency).}
If $T_i\geq k$, $Y_{ki}=Y^{(\emptyset)}_{ki}$ and $\overline{L}_{ki}=\overline{L}_{ki}^{(\emptyset)}$. $Y_i^{(T_i)}=Y_i$.
\end{assumption}

\begin{defn} \label{Hdef}
On $k\leq T_i$, define $H_i(k)=Y_{ki}$. On $k>T_i$, define $H_i(k)=Y_{ki}-\gamma_{k}^{(T_i)}\left(\overline{L}_{T_ii}\right)$.
\end{defn}
The following theorem was proven in \cite{Victor}:
\begin{theorem} {\bf (Unbiased estimating equations).}\label{see}
  Suppose that Consistency Assumption~\ref{cons} and Assumption of
  No Unmeasured Confounding~\ref{intconf} hold. With $\overline{{\cal L}}_{m}$ the space where the $\overline{L}_{mi}$ take their values, consider any
  $\vec{q}_{m}^{\;k}:\overline{{\cal L}}_{m}\rightarrow
  \mathbb{R}^{p}$, $m=0,\ldots,K$, $k>m$, which are
  measurable and bounded. Then
\begin{equation*}
E \left(\sum_{m=0}^{K}\sum_{k>m}
\vec{q}^{\;k}_{m}\left(\overline{L}_{mi}\right) H_i\left(k\right) 1_{\overline{A}_{m-1i}=\overline{0}}
\left\{A_{mi}-p\left(A_{mi}=1|\overline{L}_{mi}\overline{A}_{m-1i}\right)\right\}\right)=0,
\end{equation*}
so
\begin{equation*}
P_n\left(\sum_{m=0}^{K}\sum_{k>m} \vec{q}^{\;k}_{m}\left(\overline{L}_{m}\right) H_{\psi}\left(k\right) 1_{\overline{A}_{m-1}=\overline{0}}
\left\{A_m-p_\theta\left(A_{m}=1|\overline{L}_{m}\overline{A}_{m-1}\right)\right\}\right)=0,
\end{equation*}
stacked with
the estimating equations for $\theta$, with $P_n$ the empirical
measure $P_n(X)=1/n \sum_{i=1}^n X_i$, are unbiased estimation
equations for $(\psi,\theta)$. The $\vec{q}^{\;k}_{m}$ here are allowed to depend on $(\psi,\theta)$, as long as they are measurable and bounded for
$\left(\psi^{*},\theta^{*}\right)$.
\end{theorem}
For identifiability of $\psi$, one needs as many
estimating equations as parameters, in this case by choosing the
dimension of $\vec{q}$. A linear treatment effect model $\gamma_\psi$ such as (\ref{modelex}) leads to estimating equations that are linear in $\psi$ and thus easy to solve.

In this example, $U_{1}(Y_i,Z_i;\psi,\theta)$ equals
\begin{equation*}
\sum_{m=0}^{K}\sum_{k>m} \vec{q}^{\;k}_{m}\left(\overline{L}_{mi}\right) H_{i\psi}\left(k\right) 1_{\overline{A}_{m-1i}=\overline{0}}
\left\{A_{mi}-p_\theta\left(A_{mi}=1|\overline{L}_{mi},\overline{A}_{m-1i}\right)\right\},
\end{equation*}
and, from Theorem~\ref{see},
\begin{equation*}
E\left(U_{1}(Y_i,Z_i;\psi^*,\theta^*)\right)=0.
\end{equation*}

It can be shown that alternative estimating equations are both more efficient and doubly robust: they lead to consistent, asymptotically normal estimators when either the treatment initiation model is correctly specified or an outcome regression model is correctly specified.

\section{Regularity conditions}

We adopt some fairly usual regularity conditions.

\begin{assumption}\label{support} {\bf Regularity Condition.} We assume that $(Y,Z)$ has a density $f_{\theta,\cdot}(Y,Z)$ with respect to some measure, and its support does not depend on $\theta$.
\end{assumption}
We wrote $f_{\theta,\cdot}$ and not $f_{\psi,\theta,\cdot}$: as indicated in Section~\ref{setting}, $\psi^*$ is only defined as the solution to an equation that usually follows from identifying assumptions.

\begin{assumption}{\bf Regularity Condition.}\label{diffb} $(\psi,\theta)$ takes values in a compact space, and $U_1(\psi,\theta)$ and $U_2(\theta)$ are twice continuously differentiable in $(\psi,\theta)$.
\end{assumption}
The following regularity condition is needed to consistently estimate $\psi$ and $\theta$ based on $U_1$ and $U_2$, respectively:

\begin{assumption}{\bf Regularity Condition.}\label{inverses} $E\bigl(\left.\frac{\partial}{\partial \psi}\right|_{\psi^*} U_1(\psi,\theta^*)\bigr)$ and $E\left(\left.\frac{\partial}{\partial \theta}\right|_{\theta^*}U_2(\theta)\right)$ have an inverse.
\end{assumption}

\section{Variance estimation}

Section~\ref{Taylor} provides a variance estimator for general settings where the nuisance parameter $\theta$ is estimated based on unbiased estimating equations. Section~\ref{identnotneeded} provides an easy explanation for why estimating nuisance parameters not needed for identification of $\psi$ does not affect the variance of $\hat{\psi}$. Sections~\ref{replace}--~\ref{DR} investigate what these results imply for settings where the nuisance parameter $\theta$ is estimated by solving (partial) score equations, and where $\psi$ does not depend on $\theta$ (see Assumption~\ref{psinotftheta}). Section~\ref{small} provides a small sample correction. For completeness, Section~\ref{bootstrap} provides a proof for consistency of nonparametric bootstrap confidence intervals using Efron's percentile method, for general settings where the nuisance parameter $\theta$ is estimated based on unbiased estimating equations.

\subsection{Taylor Expansion / Middle Value Theorem and its consequences for Estimation Procedure~\ref{proc}}\label{Taylor}

The derivations in Section~\ref{Taylor} do not depend on whether or not $\theta$ is estimated solving partial score equations; $\hat{\theta}$ may solve any smooth unbiased estimating equations.

As in Chapter~5 of \cite{Vaart}, because of the Middle Value Theorem (there may be different $(\tilde{\psi},\tilde{\theta})$ between $(\psi^*,\theta^*)$ and $(\hat{\psi},\hat{\theta})$ in each row, but as usual that does not affect the limiting variance):
\begin{eqnarray*}0&=&P_n\left(\begin{array}{c} U_1(\hat{\psi},\hat{\theta})\\U_2(\hat{\theta})\end{array}\right)\\
&=&P_n\left(\begin{array}{c} U_1(\psi^*,\theta^*)\\U_2(\theta^*)\end{array}\right)+
\left(\left.\frac{\partial}{\partial(\psi,\theta)}\right|_{(\tilde{\psi},\tilde{\theta})}P_n\left(\begin{array}{c} U_1(\psi,\theta)\\U_2(\theta)\end{array}\right)\right)
\left(\begin{array}{c}\hat{\psi}-\psi^*\\\hat{\theta}-\theta^*\end{array}\right).
\end{eqnarray*}
Notice that since $(\hat{\psi},\hat{\theta})$ is consistent (Theorem~5.9 of \cite{Vaart}) and $(U_1,U_2)$ is a continuously differentiable function of $(\psi,\theta)$ (Regularity Condition~\ref{diffb}), Lemma~\ref{Donskerconv} implies that
\begin{equation*}
\left.\frac{\partial}{\partial(\psi,\theta)}\right|_{(\tilde{\psi},\tilde{\theta})}P_n\left(\begin{array}{c} U_1(\psi,\theta)\\U_2(\theta)\end{array}\right)
\rightarrow^P E\left(\left.\frac{\partial}{\partial(\psi,\theta)}\right|_{(\psi^*,\theta^*)}\left(\begin{array}{c} U_1(\psi,\theta)\\U_2(\theta)\end{array}\right)\right).
\end{equation*}
Let's study this expected derivative:
\begin{eqnarray*}E\left(\frac{\partial}{\partial(\psi,\theta)}\left(\begin{array}{c} U_1(\psi,\theta)\\U_2(\theta)\end{array}\right)\right)
&=&\left(\begin{array}{cc}E\left(\frac{\partial}{\partial \psi}U_1(\psi,\theta)\right) & E\left(\frac{\partial}{\partial \theta}U_1(\psi,\theta)\right)\\
0 & E\left(\frac{\partial}{\partial \theta}U_2(\theta)\right)\end{array}\right),
\end{eqnarray*}
so that (Partition Inverse Formula, this can easily be verified)
\begin{equation*}\left(E\frac{\partial}{\partial(\psi,\theta)}\left(\begin{array}{c} U_1(\psi,\theta)\\U_2(\theta)\end{array}\right)\right)^{-1}
=\left(\begin{array}{cc}\left(E\frac{\partial}{\partial \psi}U_1(\psi,\theta)\right)^{-1} & -B\\
0 & \left(E\frac{\partial}{\partial \theta}U_2(\theta)\right)^{-1}\end{array}\right),
\end{equation*}
with
\begin{equation*}B=\left(E\frac{\partial}{\partial \psi}U_1(\psi,\theta)\right)^{-1}\left(E\frac{\partial}{\partial \theta}U_1(\psi,\theta)\right)\left(E\frac{\partial}{\partial \theta}U_2(\theta)\right)^{-1}.
\end{equation*}
Combining, we find that $\sqrt{n}\bigl(\hat{\psi}-\psi^*\bigr)$ equals
\begin{equation*}
\left(E\frac{\partial}{\partial \psi}U_1(\psi,\theta)\right)^{-1}
\end{equation*}
\begin{equation}
\sqrt{n}\left(-P_n U_1+\left(E\frac{\partial}{\partial \theta}U_1(\psi,\theta)\right)\left(E\frac{\partial}{\partial \theta}U_2(\theta)\right)^{-1}P_nU_2\right)+o_P(1).
\label{betahatminbeta}
\end{equation}
Equation~(\ref{betahatminbeta}) leads to the setting of \cite{pierce1982asymptotic} equation~(1.2).

Equation~(\ref{betahatminbeta}) can be used to estimate the true limiting variance of $\hat{\psi}$, leading to (with $Z^{\otimes 2}=ZZ^\top$)
\begin{equation}\label{VARhatgen}
\widehat{VAR}(\hat{\psi})\approx \frac{1}{n}\left(\hat{E}\frac{\partial}{\partial \psi}U_1(\psi,\theta)\right)^{-1}
\end{equation}
\begin{equation*}
\hat{E}\left(\left(U_1-\left(E\frac{\partial}{\partial \theta}U_1(\psi,\theta)\right)\left(E\frac{\partial}{\partial \theta}U_2(\theta)\right)^{-1}U_2\right)^{\otimes 2}\right) \left(\hat{E}\frac{\partial}{\partial \psi}U_1(\psi,\theta)\right)^{-1\top}.
\end{equation*}
This variance can be estimated by plugging in $\hat{\psi}$ and $\hat{\theta}$ and replacing expectations by averages over all observations, $P_n$, to obtain 
\begin{equation}\label{VARhatgenest}
\widehat{VAR}(\hat{\psi})\approx \frac{1}{n}\left(P_n\frac{\partial}{\partial \psi}U_1(\hat{\psi},\hat{\theta})\right)^{-1}
\end{equation}
\begin{equation*}
P_n\left(\left(U_1-\left(P_n\frac{\partial}{\partial \theta}U_1(\hat{\psi},\hat{\theta})\right)\left(P_n\frac{\partial}{\partial \theta}U_2(\hat{\theta})\right)^{-1}U_2\right)^{\otimes 2}\right) \left(P_n\frac{\partial}{\partial \psi}U_1(\hat{\psi},\hat{\theta})\right)^{-1\top}.
\end{equation*}
Under regularity conditions, this leads to a consistent estimate of the limiting asymptotic variance. Implementing this variance requires programming: especially the derivative of $U_1$ with respect to $\theta$ will usually not be provided by standard software.

Ignoring estimation of the nuisance parameter $\theta$, we would have obtained
\begin{equation*}
\sqrt{n}\bigl(\hat{\psi}-\psi^*\bigr)=\left(E\frac{\partial}{\partial \psi}U_1(\psi,\theta)\right)^{-1}\sqrt{n}\left(-P_n U_1\right)+o_P(1),
\end{equation*}
leading to the (often inconsistent) approximate variance
\begin{equation}\label{VARhatincor}
\frac{1}{n}\left(\hat{E}\frac{\partial}{\partial \psi}U_1(\psi,\theta)\right)^{-1}
\hat{E}\left(U_1^{\otimes 2}\right) \left(\hat{E}\frac{\partial}{\partial \psi}U_1(\psi,\theta)\right)^{-1\top}.
\end{equation}

If $U_1$ and $U_2$ are orthogonal, a comparison of (\ref{VARhatgenest}) and (\ref{VARhatincor}) shows that the variance is typically smaller when the true $\theta$ is known and plugged in. Web-appendix~\ref{Efficient} shows this happens for efficient estimators of $\psi$.

\subsection{Estimating nuisance parameters not used for identification does not affect the variance of $\hat{\psi}$}\label{identnotneeded}

Section~\ref{identnotneeded} considers settings with additional nuisance parameters $\xi$ that appear in the estimating equations for $\psi$ and that, unlike $\theta$, are not needed for identifiability; that is,
$EU_1(\psi^*,\theta^*,\xi)=0$ for all $\xi$. For example, $\xi$ could be included to increase precision.

Under regularity conditions, we then have \begin{equation*}E\frac{\partial}{\partial \xi} U_1(\psi^*,\theta^*,\xi)=\frac{\partial}{\partial \xi}EU_1(\psi^*,\theta^*,\xi)=0.
\end{equation*}
Suppose furthermore that also $\xi$ is estimated in a first step, and solves unbiased estimating equations $P_nU_3(\xi)=0$ with $EU_3(\xi^*)=0$.
In this setting, the second term in equation~(\ref{betahatminbeta}) becomes
\begin{equation*}
\sqrt{n}\left(\begin{array}{cc}E\frac{\partial}{\partial \theta}U_1(\psi,\theta) & 0\end{array}\right)\left(\begin{array}{cc}E\frac{\partial}{\partial \theta}U_2(\theta)&0\\ 0& E\frac{\partial}{\partial \xi}U_3(\xi)\end{array}\right)^{-1}\left(\begin{array}{c}P_nU_2\\ P_nU_3\end{array}\right),
\end{equation*}
or equivalently
\begin{eqnarray}
\lefteqn{\sqrt{n}\left(\begin{array}{cc}E\frac{\partial}{\partial \theta}U_1(\psi,\theta) & 0\end{array}\right)\left(\begin{array}{c}\left(E\frac{\partial}{\partial \theta}U_2(\theta)\right)^{-1}P_nU_2\\ \left(E\frac{\partial}{\partial \xi}U_3(\xi)\right)^{-1}P_nU_3\end{array}\right)}\nonumber\\
&=&\sqrt{n}\left(E\frac{\partial}{\partial \theta}U_1(\psi,\theta)\right) \left(E\frac{\partial}{\partial \theta}U_2(\theta)\right)^{-1}P_nU_2
.
\label{betahatminbetaxi}\end{eqnarray}
Thus, all contributions of $\xi$ disappear, and equation~(\ref{betahatminbeta}) and all its consequences remain valid: the asymptotic variance of $\hat{\psi}$ does not depend on whether $\xi$ is estimated or known.

\subsection{Replacing the estimating equations without changing the estimator $(\hat{\psi},\hat{\theta})$}\label{replace}

Observation: we actually solve the equations not only $P_n(U_1, U_2)=0$, but then also automatically all linear combinations $P_n(U_1+CU_2)=0$ with matrix $C$ fixed. In particular, if $A_{\psi^*,\theta^*}U_2$ is the projection of $U_1$ onto the space spanned by $U_2$, we equivalently to (\ref{stack}) solve
\begin{equation*}
\left(\begin{array}{c}P_n \left(U_1(\psi,\theta)-A_{\psi^*,\theta^*}U_2(\theta)\right)\\P_nU_2(\theta)\end{array}\right)=0.
\end{equation*}
Notice that we use $A_{\psi^*,\theta^*}$ here, and not $A_{\psi,\theta}$; if we would use $A_{\psi,\theta}$, we would change the estimator $\hat{\psi}$.

It is easy to calculate $A_{\psi^*,\theta^*}$:
\begin{eqnarray*}
0&=&E\left(\left(U_1(\psi,\theta)-A_{\psi,\theta}U_2(\theta)\right)U_2^\top\right)\\
&=&E\left(\left(U_1(\psi,\theta)U_2(\psi,\theta)^\top-A_{\psi,\theta}U^{\otimes 2}_2(\theta)\right)\right),
\end{eqnarray*}
so
\begin{equation}\label{Apsitheta}
A_{\psi,\theta}=E\left(U_1(\psi,\theta)U_2(\psi,\theta)^\top\right)\left(E\left(U^{\otimes 2}_2(\theta)\right)\right)^{-1}.
\end{equation}

Now consider the projection $\tilde{U}_1$ of $U_1$ onto the orthocomplement of $U_2$:
\begin{eqnarray}\label{U1tilde}
\lefteqn{\tilde{U}_1(\psi,\theta)=U_1(\psi,\theta)-A_{\psi^*,\theta^*}U_2(\theta)}\nonumber\\
&=&U_1-E\left(U_1(\psi^*,\theta^*)U_2(\psi^*,\theta^*)^\top\right)\left(E\left(U^{\otimes 2}_2(\theta^*)\right)\right)^{-1}U_2;
\end{eqnarray}
this is just another $U_1$, leading to the same $\hat{\psi}$. Thus, like in the derivations that led to equation~(\ref{betahatminbeta}), we find that $\sqrt{n}\bigl(\hat{\psi}-\psi^*\bigr)$ equals
\begin{equation*}
\left(E\frac{\partial}{\partial \psi}\tilde{U}_1(\psi,\theta)\right)^{-1}
\end{equation*}
\begin{equation}
\sqrt{n}\left(-P_n \tilde{U}_1+\left(E\frac{\partial}{\partial \theta}\tilde{U}_1(\psi,\theta)\right)\left(E\frac{\partial}{\partial \theta}U_2(\theta)\right)^{-1}P_nU_2\right)+o_P(1).\label{betahatminbetatildetilde}
\end{equation}

Notice that $\tilde{U}_1$ depends on $U_1$ and a linear combination of $U_2$ that depends on the true $\psi^*$, not on $\psi$. Thus, the derivative of $\tilde{U}_1$ with respect to $\psi$ is the same as the derivative of $U_1$ with respect to $\psi$, and (\ref{betahatminbetatildetilde}) implies that $\sqrt{n}\bigl(\hat{\psi}-\psi^*\bigr)$ equals
\begin{equation*}
\left(E\frac{\partial}{\partial \psi}U_1(\psi,\theta)\right)^{-1}
\end{equation*}
\begin{equation}
\sqrt{n}\left(-P_n \tilde{U}_1+\left(E\frac{\partial}{\partial \theta}\tilde{U}_1(\psi,\theta)\right)\left(E\frac{\partial}{\partial \theta}U_2(\theta)\right)^{-1}P_nU_2\right)+o_P(1).
\label{betahatminbetatilde}\end{equation}
Equation~(\ref{betahatminbetatilde}) does not depend on $\hat{\theta}$ solving (partial) score equations.

\subsection{$\theta$ estimated by partial score equations, $\psi$ does not depend on $\theta$: consequences of equivalent estimating equations}\label{consequences}

For $\hat{\theta}$ solving (partial) score equations (Assumption~\ref{MPLE}) and if Assumption~\ref{psinotftheta} holds, a funny thing happens. In such setting,
\begin{equation}\label{derivative0}
E_{\psi,\theta}\left(\frac{\partial}{\partial \theta}\tilde{U}_1(\psi,\theta)\right)=0.
\end{equation}
To see this, first recall that since both $U_1$ and $U_2$ lead to unbiased estimating equations,
\begin{equation*}
E_{\psi,\theta}\left(\tilde{U}_1(\psi,\theta)\right)=0.
\end{equation*}
To take the derivative of this expression with respect to $\theta$ while keeping $\psi$ fixed, we need Assumption~\ref{psinotftheta}; if Assumption~\ref{psinotftheta} does not hold the derivative with respect to $\theta$ of this expression would also involve $\frac{\partial}{\partial \theta}\psi(\theta)$. Under Assumption~\ref{psinotftheta},
\begin{eqnarray*}
0&=&\frac{\partial}{\partial \theta}E_{\theta}\left(\tilde{U}_1(\psi,\theta)\right)
=\frac{\partial}{\partial \theta}\int \tilde{U}_1(\psi,\theta) f_{\theta,\cdot}(y,z) d(y,z)\\
&=&\int \left(\frac{\partial}{\partial \theta}\tilde{U}_1(\psi,\theta)\right) f_{\theta,\cdot}(y,z) d(y,z)+\int \tilde{U}_1(\psi,\theta) \left(\frac{\partial}{\partial \theta}f_{\theta,\cdot}(y,z)\right) d(y,z)\\
&=&E\left(\frac{\partial}{\partial \theta}\tilde{U}_1(\psi,\theta)\right) + \int \tilde{U}_1(\psi,\theta) \left(\frac{\partial}{\partial \theta}\log f_{\theta,\cdot}(y,z)\right)  f_{\theta,\cdot}(y,z) d(y,z)\\
&=&E\left(\frac{\partial}{\partial \theta}\tilde{U}_1(\psi,\theta)\right)+E\left(\tilde{U}_1(\psi,\theta)U_2^\top(\theta)\right)
=E\left(\frac{\partial}{\partial \theta}\tilde{U}_1(\psi,\theta)\right).
\end{eqnarray*}
For the third equality we use that the support of $(Y,Z)$ does not depend on $\theta$ (Regularity Condition~\ref{support}), and for the last equality we use that $\tilde{U}_1$ and $U_2$ are orthogonal. Equation~(\ref{derivative0}) follows.

Combining equations~(\ref{derivative0}) and~(\ref{betahatminbetatilde}), it follows that
\begin{equation}\label{betahatminbetaprime}
\sqrt{n}\bigl(\hat{\psi}-\psi^*\bigr)
=-E\left(\frac{\partial}{\partial \psi}U_1(\psi,\theta)\right)^{-1}\sqrt{n} P_n \tilde{U}_1+o_P(1),
\end{equation}
with (by construction) $\tilde{U}_1$ the projection of $U_1$ onto the orthocomplement of $U_2$.

\subsection{True limiting variance of $\hat{\psi}$ is smaller if $\theta$ is estimated using (partial) score equations, $\psi$ does not depend on $\theta$}\label{nuisestVAR}

Let's contrast equation~(\ref{betahatminbetaprime}) with
a setting where we would have (somehow) plugged in the true nuisance parameter $\theta$ into $U_1$, leading to an estimator $\tilde{\psi}$ with (using the reasoning from Section~\ref{Taylor} or \cite{Vaart} Theorem~5.21)
\begin{equation*}
\sqrt{n}\left(\tilde{\psi}-\psi^*\right)=-\left(E\frac{\partial}{\partial \psi}U_1(\psi,\theta^*)\right)^{-1}\sqrt{n}P_n U_1+o_P(1).
\end{equation*}
From equation~(\ref{betahatminbetaprime}), since $\tilde{U}_1$ from Sections~\ref{replace} and~\ref{consequences} is a projection of $U_1$, the limiting variance of $\hat{\psi}$ is \emph{smaller} when the nuisance parameter $\theta$ is estimated by solving (partial) score equations!

To be precise: the limiting asymptotic variance when the nuisance parameters $\theta$ are not estimated is
\begin{equation}\label{notest}
E\left(\frac{\partial}{\partial \psi}U_1(\psi,\theta)\right)^{-1}EU_1^{\otimes 2}\left(\psi,\theta\right)E\left(\frac{\partial}{\partial \psi}U_1(\psi,\theta)\right)^{-1\top}.
\end{equation}
The variance of $\hat{\psi}$ when $\theta$ is estimated by solving (partial) score equations is, from equation~(\ref{betahatminbetaprime}),
\begin{equation}\label{est}
E\left(\frac{\partial}{\partial \psi}U_1(\psi,\theta)\right)^{-1}E\tilde{U}_1^{\otimes 2}\left(\psi,\theta\right)E\left(\frac{\partial}{\partial \psi}U_1(\psi,\theta)\right)^{-1\top}.
\end{equation}
The latter is the smaller variance, since $\tilde{U}_1$ is a projection of $U_1$, so that $EU_1^{\otimes 2}-E\tilde{U}_1^{\otimes 2}$ is positive semi-definite.

\subsection{Estimating $\theta$ using (partial) score equations and ignoring that, leads to conservative inference for $\hat{\psi}$ if $\psi$ does not depend on $\theta$}\label{conservative}

Often, one will first estimate $\theta$ with standard software (solving (partial) score equations), then plug this $\hat{\theta}$ into standard software to estimate $\psi$. This leads to $\hat{\psi}$, see Equation~(\ref{EE}). One could be tempted to then use the confidence intervals for $\psi$ provided by the standard software, as in for example \cite{MSM1} and \cite{MSM2}. If the sandwich estimator of the variance is used to construct the confidence intervals for $\psi$, this leads to conservative inference, since the sandwich estimator will use
\begin{equation}\label{sandwichest}
\hat{E}\left(\frac{\partial}{\partial \psi}U_1(\hat{\psi},\hat{\theta})\right)^{-1}\hat{E}U_1^{\otimes 2}\left(\hat{\psi},\hat{\theta}\right)\hat{E}\left(\frac{\partial}{\partial \psi}U_1(\hat{\psi},\hat{\theta})\right)^{-1\top},
\end{equation}
which is consistent for the variance in the setting where we plug in the true $\theta$ (see equation~(\ref{notest})). Section~\ref{nuisestVAR} shows that this is a larger variance than the variance (\ref{betahatminbetaprime}) when $\theta$ is estimated: conservative inference.

\subsection{Consistent estimates for the variance of $\hat{\psi}$ if $\hat{\theta}$ solves (partial) score equations, $\psi$ does not depend on $\theta$}

Starting with equation~(\ref{est}), we estimate the asymptotic variance consistently 
by plugging in $\hat{\psi}$ and $\hat{\theta}$ and replacing expectations by empirical averages over all observations, $P_n$, to obtain 
\begin{equation}\label{ASVAR}
P_n\left(\frac{\partial}{\partial \psi}U_1(\hat{\psi},\hat{\theta})\right)^{-1}P_n\tilde{U}_1^{\otimes 2}\left(\hat{\psi},\hat{\theta}\right)P_n\left(\frac{\partial}{\partial \psi}U_1(\hat{\psi},\hat{\theta})\right)^{-1\top}.
\end{equation}
Under regularity conditions, this is a consistent estimate of the limiting asymptotic variance.

Expression~(\ref{ASVAR}) involves $\tilde{U}_1$, which we then need to estimate. Recall from equation~(\ref{U1tilde}) that
\begin{equation*}
\tilde{U}_1(\psi,\theta)=
U_1-E\left(U_1(\psi^*,\theta^*)U_2(\theta^*)^\top\right)\left(E\left(U^{\otimes 2}_2(\theta^*)\right)\right)^{-1}U_2.
\end{equation*}
$U_1$ can be plugged in, and $E\left(U^{\otimes 2}_2(\theta^*)\right)$ is the (partial) Fisher information for $\theta$, to be obtained from estimating $\theta$. Depending on computing considerations, one can choose to either estimate
$E\left(U_1(\psi^*,\theta^*)U_2(\theta^*)^\top\right)$ or $E\frac{\partial}{\partial \theta}U_1(\psi^*,\theta^*)$, because as shown in Web-appendix~A,
\begin{equation}
E\frac{\partial}{\partial \theta}U_1(\psi,\theta)
=-EU_1U_2^\top.\label{ddtheta}
\end{equation}
Of these, $\hat{E}\left(U_1(\psi^*,\theta^*)U_2(\theta^*)^\top\right)$ may often be easier to program: one could plug in the estimating equations for $\theta$ and $\psi$, replace $(\psi^*,\theta^*)$ by $(\hat{\psi},\hat{\theta})$, and replace $E$ by the empirical average $P_n$.

\subsection{Estimating VAR($\hat{\psi}$) when $\hat{\theta}$  solves (partial) score equations, $\psi$ does not depend on $\theta$: building on standard software}\label{varcorrection}

The variance of $\hat{\psi}$ can also be estimated by subtracting a variance correction term from the usual sandwich estimator, often provided by standard software after the nuisance parameter $\theta$ is estimated in a first step. This works as follows.
In Web-appendix~A we show that equation~(\ref{betahatminbetaprime}) implies that the variance of $\hat{\psi}$ is approximately
\begin{eqnarray}\label{g-var}
\lefteqn{\frac{1}{n}E\left(\frac{\partial}{\partial \psi}U_{1}(\psi,\theta)\right)^{-1}\left(E\left(U_{1}^{\otimes 2}\right)
-E\left(U_{1}U_{2}^\top\right)
E\left(U_{2}^{\otimes 2}\right)^{-1}E\left(U_{1}U_{2}^\top\right)^\top\right)}\nonumber\\
&&\;\hspace{5cm}\left(E\left(\frac{\partial}{\partial \psi}U_{1}(\psi,\theta)\right)\right)^{-1\,\top}
\end{eqnarray}
in larger samples. The first term of this variance of $\hat{\psi}$, 
\begin{equation}
\frac{1}{n}E\left(\frac{\partial}{\partial \psi}U_{1}(\psi,\theta)\right)^{-1}\left(EU_{1}^{\otimes 2}(\psi,\theta)\right) E\left(\frac{\partial}{\partial \psi}U_{1}(\psi,\theta)\right)^{-1\,\top},\label{sandwich}
\end{equation}
is the result that most software will provide as the ``sandwich estimator'' or ``robust estimator'' for the variance of $\hat{\psi}$ when the estimated $\hat{\theta}$ is provided to solve equation~(\ref{EE}). 

Thus, we can focus on the variance correction term in equation~(\ref{g-var}),
\begin{eqnarray}\label{g-asvar}
\lefteqn{-\frac{1}{n}E\left(\frac{\partial}{\partial \psi}U_{1}(\psi,\theta)\right)^{-1}
E\left(U_{1}U_{2}^\top\right)
E\left(U_{2}^{\otimes 2}\right)^{-1}E\left(U_{1}U_{2}^\top\right)^\top}\nonumber\\
&&\;\hspace{5cm}\left(E\left(\frac{\partial}{\partial \psi}U_{1}(\psi,\theta)\right)\right)^{-1\,\top}.
\end{eqnarray}
An estimate of $E\left(U_{2}^{\otimes 2}\right)$ can often be obtained from standard software used to estimate $\theta$, using that standard software usually estimates
\begin{equation*}
\widehat{{\rm VAR}}(\hat{\theta})=\frac{1}{n}\left(\hat{E}\left(U_{2}^{\otimes 2}\right)\right)^{-1},
\end{equation*}
based on classical theory on the (partial) Fisher information for $\theta$.

Thus, we can estimate the variance of $\hat{\psi}$ by the sandwich estimator of the variance ignoring estimation of $\theta$ minus
\begin{eqnarray}\label{varcorterm}
\lefteqn{\hat{E}\left(\frac{\partial}{\partial \psi}U_{1}(\psi,\theta)\right)^{-1}
\hat{E}\left(U_{1}U_{2}^\top\right)
\widehat{{\rm VAR}}(\hat{\theta}) \hat{E}\left(U_{1}U_{2}^\top\right)^\top}\nonumber\\
&&\;\hspace{5cm} \left(\hat{E}\left(\frac{\partial}{\partial \psi}U_{1}(\psi,\theta)\right)\right)^{-1\,\top}.
\end{eqnarray}

In some cases, standard software estimating $\psi$ may provide the ``bread'' $\hat{E}\left(\frac{\partial}{\partial \psi}U_{1}(\psi,\theta)\right)^{-1}$ of this sandwich estimator correction, which is the same as in equation~(\ref{sandwich}); in most cases, $E\left(U_{1}U_{2}^\top\right)$ will have to be estimated by for example $P_n U_1 U_2^\top$ with $U_1$ and $U_2$ estimated by plugging in $(\hat{\psi},\hat{\theta})$.

\subsection{Is it truly conservative? Can the variance correction term be zero when $\theta$ is estimated using (partial) score equations, $\psi$ does not depend on $\theta$?}


If $\theta$ is estimated using (partial) score equations, the (conservative) sandwich estimator for the variance of $\hat{\psi}$ ignoring estimation of $\theta$, that is, ignoring the variance correction term (\ref{varcorterm}), leads to a consistent estimator for the variance if and only if the estimating equations for $\psi$ based on $U_1$ are orthogonal to the (partial) score equations for $\theta$. That is, if and only if $U_1$ and $U_2$ are orthogonal: $EU_1U_2^\top=0$. 

This can be seen as follows. First, by equation~(\ref{Apsitheta}) or general theory, the projection of $U_1$ onto the space spanned by $U_2$ equals
\begin{equation*} E(U_1U_2^\top) (EU_2^{\otimes 2})^{-1} U_2.
\end{equation*}
Second, the variance of this projection equals
\begin{eqnarray*}
\lefteqn{E(U_1U_2^\top) (EU_2^{\otimes 2})^{-1}EU_2^{\otimes 2} (EU_2^{\otimes 2})^{-1}E(U_1U_2^\top)^\top}\\
&=&E(U_1U_2^\top) (EU_2^{\otimes 2})^{-1}E(U_1U_2^\top)^\top.
\end{eqnarray*}
Since this is the variance of the projection, it equals $0$ if and only if that projection is $0$ with probability $1$, that is, if $U_1$ and $U_2$ are orthogonal. But this variance is the inside of the variance correction term (\ref{varcorterm}), and to obtain (\ref{varcorterm}), we pre- and post-multiply by an invertible matrix (Regularity Condition~\ref{inverses}). Thus, the variance correction term equals $0$ if and only if $U_1$ and $U_2$ are orthogonal.

\subsection{Estimating $\theta$ using (partial) score equations does not affect the variance of efficient estimators for $\psi$ if $\psi$ does not depend on $\theta$}\label{efficient}

If an estimator for $\psi$ is used that would be efficient were $\theta$ known, solving unbiased estimating equations of the form $P_n U^{eff}_1(\psi,\theta^*)=0$, the variance of $\hat{\psi}$ is not affected by whether the nuisance parameter $\theta$ is estimated by solving (partial) score equations or the true $\theta$ is plugged in. That means that for such estimator, the true limiting variance is given by the sandwich estimator (\ref{sandwichest}), which will typically be provided by standard software if $\hat{\theta}$ is plugged in as a first step.


Providing more detail, suppose we have an efficient estimator $\hat{\tilde{\psi}}$ for $\psi$ when $\theta$ is known, solving unbiased estimating equations $P_n U^{eff}_1(\psi,\theta^*)=0$. If $U^{eff}_1(\psi,\hat{\theta})$ is used to estimate $\psi$, with $\hat{\theta}$ solving (partial) score equations for $\theta$, then 1. the variance of $\hat{\psi}$ is the same regardless of whether $\hat{\theta}$ is estimated or the true $\theta$ is plugged in and 2. the sandwich estimator for the variance of $\hat{\psi}$ based on $U_1$ with $\hat{\theta}$ plugged in as a first step is consistent.

The idea behind the proof of this result is to create more efficient estimating equations by projecting $U_1^{eff}$ onto the orthocomplement of the space spanned by $U_2$. Since $U_1^{eff}$ is efficient, this should not lead to an efficiency gain. Thus, $U_1^{eff}$ is already orthogonal to $U_2$, and so there is no efficiency gain from estimating $\theta$. Carrying this out directly is complicated, since the projection of $U_1^{eff}$ onto the orthocomplement of the space spanned by $U_2$ depends on the true parameter $\psi^*$ and two expectations. Web-appendix~B shows this can be solved, by estimating those in a first step.

\subsection{Variance estimation of doubly robust estimators}
\label{DR}

Our results also have implications for variance estimation of doubly robust estimators. For example, the estimators for the effect of a binary treatment in \cite{Bang} are doubly robust: they are consistent and asymptotically normal if either the propensity score model or an outcome regression model is correctly specified. If the outcome regression model is misspecified but the treatment probabilities are correctly specified, the sandwich estimator for the variance ignoring estimation of $\theta$ is conservative. If both models are correctly specified this sandwich estimator of the variance is consistent: since neither model is needed separately to identify $\psi$, it leads to the setting in Section~\ref{identnotneeded}. 

This is a general pattern for doubly robust estimators, where if either the model for nuisance parameter $\theta_1$ or the model for nuisance parameter $\theta_2$ is correctly specified, the estimating equations for $\psi$ are unbiased. In such settings, if one of the nuisance parameter models is correctly specified and estimated by solving (partial) score equations and $\psi$ does not depend on that parameter, the sandwich estimator is conservative. If both models are correctly specified, the sandwich estimator is consistent, since neither model is needed separately to identify $\psi$, leading to the setting in Section~\ref{identnotneeded}.

\subsection{Small sample corrections}\label{small}

It is well-known that the sandwich estimator of the variance may underestimate the true variance (e.g., \cite{li2015small}, \cite{mackinnon1985some}). Confidence interval coverage can be improved by multiplying the sandwich estimator of the variance by $n/(n-\dim \psi)$ or $n/(n-\dim \psi-\dim \theta)$ (df- or degrees-of-freedom correction, \cite{li2015small}, \cite{mackinnon1985some}) and using quantiles of the t-distribution with degrees of freedom equal to $n$ minus the dimension of the parameters involved instead of quantiles of the normal distribution (\cite{li2015small}). We implemented this in the application sections. Simulations in \cite{li2015small} and \cite{mackinnon1985some} show that other methods may work better in specific applications, depending on the form of $U_1$. Developing such methods in settings with nuisance parameters, probably depending on the form of $U_1$ and possibly also $U_2$, is an interesting topic for future research.

\subsection{The bootstrap for consistent variance estimation}\label{bootstrap}

Confidence intervals for $\psi$ could also be constructed using the nonparametric bootstrap with Efron's percentile method (\cite{Vaart}).
The nonparametric bootstrap often leads to confidence intervals with asymptotically correct coverage because of the Bootstrap Master Theorem in \cite{Kosorokbook}: the class of stacked estimating functions is Donsker because of Example~19.7 in \cite{Vaart} and Regularity Conditions~\ref{diffb}, provided that for both $j=1$ and $j=2$,
\begin{equation*}\left|U_j(Y,Z;\psi_1, \theta_1)-U_j(Y,Z;\psi_2, \theta_2)\right|\leq m_j(Y,Z)\left\|(\theta_1,\psi_1)-(\theta_2,\psi_2)\right\| 
\end{equation*}
for every $\psi_1, \theta_1,\psi_2,\theta_2$ for some $m_j$ with $E\bigl(m_j^2\bigr)<\infty$. This condition is satisfied for example if the derivative of $U_1$ with respect to $(\psi,\theta)$ and the derivative of $U_2$ with respect to $\theta$ are bounded.

This implies that the bootstrap indeed picks up estimation of the nuisance parameter $\theta$, intuitively since the bootstrap estimates $\theta$ in each bootstrap sample. Unfortunately, with complicated estimation procedures and/or large datasets the bootstrap can take considerable computing time, as in for example Section~\ref{ART}.

%


\section{The effect of cazavi versus colistin on resistant bacterial infections: Inverse Probability of Treatment Weighting for point treatment}\label{cazavi}

Resistant bacterial infections are an increasing threat to public health, with limited new antibiotics reaching market (\cite{renwick2016systematic}). 
\cite{Duin} estimated the effect of a relatively new antibiotic, cazavi (ceftazidime-avibactam), versus colistin on resistant bacterial infections using Inverse Probability of Treatment Weighting (IPTW). Confidence intervals were estimated using the semiparametric bootstrap and Efron's percentile method, fixing the number of patients in each treatment group as proposed in for example \cite{imbens2015causal}. The analyses included $n=99$ patients treated with colistin and  $n=38$ patients treated with cazavi. The primary efficacy outcome was categorical: hospital death in the first 30 days/discharged home in the first 30 days/other.

As described in Section~\ref{IPTWpoint}, Inverse Probability of Treatment Weighting (IPTW) for point treatment has $U_2$ as in equation~(\ref{IPTWpointU2}) and
\begin{eqnarray*}
U_1(\psi_1,\psi_0,\hat{\theta})&=&\left(\begin{array}{c} \frac{1_{A_i=1}}{P_{\hat{\theta}}\left(A_i=1|L_i\right)}\left(Y_i-\psi_1\right)\\
\frac{1_{A_i=0}}{P_{\hat{\theta}}\left(A_i=0|L_i\right)}\left(Y_i-\psi_0\right)\end{array}\right),
\end{eqnarray*}
leading to 
\begin{equation*}
\hat{\psi}_1=\frac{P_n Y \frac{1_{A=1}}{P_{\hat{\theta}}\left(A=1|L\right)}}{P_n \frac{1_{A=1}}{P_{\hat{\theta}}\left(A=1|L\right)}},\;\;\;\;\;\;\;\hat{\psi}_0=\frac{P_n Y \frac{1_{A=0}}{P_{\hat{\theta}}\left(A=0|L\right)}}{P_n \frac{1_{A=0}}{P_{\hat{\theta}}\left(A=0|L\right)}},
\end{equation*}
or $(\hat{\psi}_1,\hat{\psi}_0)$ solving the estimating equations
\begin{equation*}C_1\left(\begin{array}{c}\psi_1\\ \psi_0 \end{array}\right)-C_2=0
\end{equation*}
with
\begin{equation*}
C_1=\left(\begin{array}{cc}P_n\frac{1_{A=1}}{P_{\hat{\theta}}\left(A=1|L\right)} & 0\\ 0& P_n\frac{1_{A=0}}{P_{\hat{\theta}}\left(A=0|L\right)}\end{array}\right),\;\;\;\;\;\;
C_2=\left(\begin{array}{c}P_n\frac{1_{A=1}}{P_{\hat{\theta}}\left(A=1|L\right)}Y\\ P_n\frac{1_{A=0}}{P_{\hat{\theta}}\left(A=0|L\right)}Y\end{array}\right).
\end{equation*}
Thus, the conservative sandwich estimator of the variance of $(\hat{\psi}_1,\hat{\psi}_0)$ ignoring estimation of $\theta$ as in Equation~(\ref{sandwich}) equals \begin{equation*}\label{sandwichIPTW}\frac{1}{n}C_1^{-1}\hat{D}C_1^{-1\top}
\end{equation*}
with
$D=EU_1^{\otimes 2}(\psi^*,\theta^*)$ which is typically estimated by $P_n U_1^{\otimes 2}(\hat{\psi},\hat{\theta})$. In this case, $D$ can be simplified, since the two components $U_{11}$ and $U_{10}$ of $U_1$ are orthogonal for the true parameters $(\psi,\theta)$: they each have expectation $0$ and
\begin{equation*}
E\left(\frac{1_{A_i=1}}{P_{\hat{\theta}}\left(A_i=1|L_i\right)}\left(Y_i-\psi_1\right)
\frac{1_{A_i=0}}{P_{\hat{\theta}}\left(A_i=0|L_i\right)}\left(Y_i-\psi_0\right)\right)=E(0)=0.
\end{equation*} 
Therefore, $D$ is orthogonal, and the resulting conservative variance estimator is diagonal with 
\begin{equation}
    \widehat{VAR}(\hat{\psi}_1)=\frac{1}{n}\left(P_n\frac{1_{A_i=1}}{P_{\hat{\theta}}\left(A_i=1|L_i\right)}\right)^{-2}P_n\left(\frac{1_{A_i=1}}{P^2_{\hat{\theta}}\left(A_i=1|L_i\right)}\left(Y_i-\hat{\psi}_1\right)^2\right),
\end{equation}
with a similar expression for $\widehat{VAR}(\hat{\psi}_0)$. With $w_{i1}=1_{A_i=1}/P_{\hat{\theta}}\left(A_i=1|L_i\right)$ and $w_{i0}=1_{A_i=0}/P_{\hat{\theta}}\left(A_i=0|L_i\right)$, this becomes
\begin{equation}
    \widehat{VAR}(\hat{\psi}_1)=\frac{1}{n}\left(\frac{\sum_{i=1}^n w_{i1}}{n}\right)^{-2}\left(\frac{\sum_{i=1}^n w_{i1}^2 (Y_i-\hat{\psi}_1)^2}{n}\right)=\frac{\sum_{i=1}^n w_{i1}^2 (Y_i-\hat{\psi}_1)^2}{(\sum_{i=1}^n w_{i1})^2},
\end{equation}
with a similar expression for $\widehat{VAR}(\hat{\psi}_0)$.
From Section~\ref{conservative}, this leads to a conservative estimator of the variances of $\hat{\psi}_1$ and $\hat{\psi}_0$.

The variance correction term (\ref{varcorterm}) of Section~\ref{varcorrection} equals
\begin{equation*}
    C_1^{-1} (P_n U_1 U_2^\top) \widehat{VAR}(\hat{\theta}) (P_n U_1 U_2^\top)^{\top} C_1^{-1},
\end{equation*}
with 
\begin{equation*}
    C_1^{-1} =\left(\begin{array}{cc} (\sum_{i=1}^n w_{i1}/n)^{-1}& 0\\ 0& (\sum_{i=1}^n w_{i0}/n)^{-1}\end{array}\right).
\end{equation*}

If one wants to avoid that the confidence intervals for the probabilities $\psi_1$ and $\psi_0$ may include $0$ or $1$, one could as usual create confidence intervals for the log oddses first and transform those to confidence intervals for the probabilities. With the log oddses $\beta_1=\log(\psi_1/(1-\psi_1))$ and $\beta_0=\log(\psi_0/(1-\psi_0))$, this leads to
\begin{equation*}
    U_1(Y_i,Z_i;\beta,\theta)=\left(\begin{array}{c}w_{1i}\left(Y_i-\frac{e^{\beta_1}}{1+e^{\beta_1}}\right)\\
    w_{0i}\left(Y_i-\frac{e^{\beta_0}}{1+e^{\beta_0}}\right)\end{array}\right).
\end{equation*}
From straightforward calculations, the conservative sandwich estimator for the variance of $\hat{\beta}$ ignoring estimation of $\theta$ equals
\begin{equation*}
\widehat{VAR}(\beta)=  C_3^{-1}  \hat{D} C_3^{-1}/n
\end{equation*}
with $\hat{D}=P_nU_1^{\otimes 2}(\hat{\beta},\hat{\theta})$ as before and
\begin{equation*}
    C_3 =-\left(\begin{array}{cc} \hat{\psi}_1(1-\hat{\psi}_1)(P_n w_{1})& 0\\ 0& \hat{\psi}_0(1-\hat{\psi}_0)(P_n w_{0})\end{array}\right).
\end{equation*}
The variance correction term (\ref{varcorterm}) of Section~\ref{varcorrection} for $\beta$ equals
\begin{equation*}
    C_3^{-1} (P_n U_1 U_2^\top) \widehat{VAR}(\hat{\theta}) (P_n U_1 U_2^\top)^{\top} C_3^{-1}.
\end{equation*}

Table~1 includes the results for the probability of hospital death in the first 30 days, the probability of being discharged home in the first 30 days, and the probability of being alive in the hospital or discharged not-to-home in the first 30 days, under cazavi and colistin. The model for treatment was a logistic regression model with predictors Pitt score less than 4 (yes/no) and infection type (Bloodstream/Urinary tract/Other).

In the Crackle-1 data, the weights ranged $1.2-5.2$. 
Table~1 reports the nonparametric bootstrap, because of Section~\ref{bootstrap} (with weights truncated as is often done in practice (\cite{Cole}), at 15, assuming all patients had a probability of at least 1 in 15 of being given each of the treatments). 

The differences between the conservative, small sample corrected consistent, and consistent confidence intervals were small. For the largest treatment group, colistin ($n=99$), the bootstrap confidence intervals were very close to the sandwich-based confidence intervals. The smallest treatment group, cazavi ($n=38$) turns out too small to observe such similarities, which is also reflected in some of the confidence intervals for the difference between colistin and cazavi.

\begin{center}
   \includegraphics[scale=1.0,angle=0]{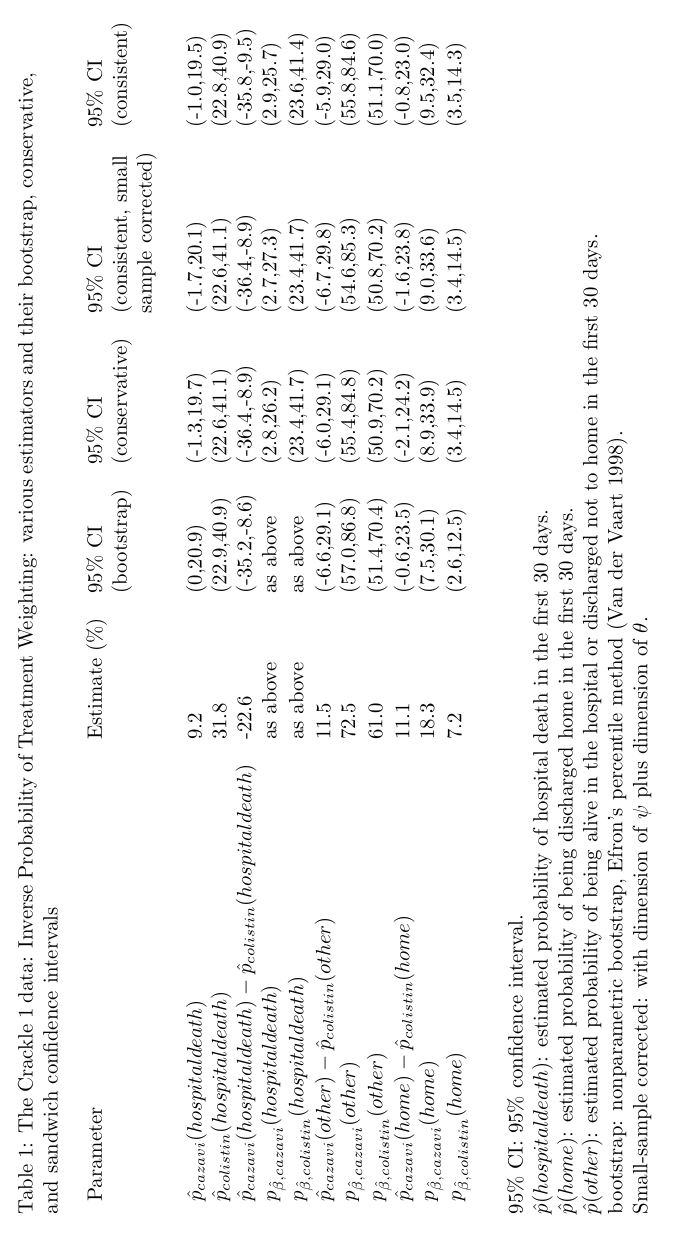}
\end{center}


\section{The effect of caz-avi versus colistin on resistant bacterial infections: simulations}\label{simulations}

We simulated 8,000 datasets with 130 and 1000 patients each, modeled after the application of Section~\ref{cazavi}. We simulated 2 independent baseline covariates: Pitt score less than 4 in 43\% of patients, and infection type: Bloodstream in 46\%, Urinary tract in 14\%, and Other in 40\% of patients. The simulated probabilities of hospital death under cazavi and colistin were
\begin{equation*}
p_{cazavi}(hospitaldeath)=\frac{1}{1+exp(1.4+11*Pittless4+0.56*bloodstream+0.28*urinary)},    
\end{equation*}
\begin{equation*}
p_{colistin}(hospitaldeath)=\frac{1}{1+exp(0.20+2.0*Pittless4-0.32*bloodstream+0.89*urinary)}.    
\end{equation*}
In the simulations, the probability of receiving caz-avi (versus colistin) was
\begin{equation*}
P(cazavi)=\frac{1}{1+exp(1-0.62*Pittless4+0.44*bloodstream+0.33*urinary)}.    
\end{equation*}
The true simulated underlying probability of hospital death was 9.0\% (0.0900182)
under caz-avi (based on 10000 times 1000 hospital death probabilities simulated under caz-avi generated with a different seed), and 30.8\% (0.3076207)
under colistin (based on 10000 times 1000 hospital death probabilities simulated under colistin).

These 8000 datasets were analyzed using Inverse Probability of Treatment Weighting (IPTW), with the model for treatment a logistic regression model with predictors Pitt score less than 4 (yes/no) and infection type (Bloodstream/Urinary tract/Other). For the setting with n=130 patients in each dataset, some treatment models resulted in a quasi-complete separation of data points. 6.3\% of estimates for the probability of hospital death under caz-avi were 0 since those datasets had no hospital deaths in this treatment group; for those datasets the true probability was not in the conservative 95\% confidence interval and the consistent estimator of the variance turned out negative. It is likely though that for such data, a different analysis strategy would be used.

\begin{center}
   \includegraphics[scale=1.0,angle=0]{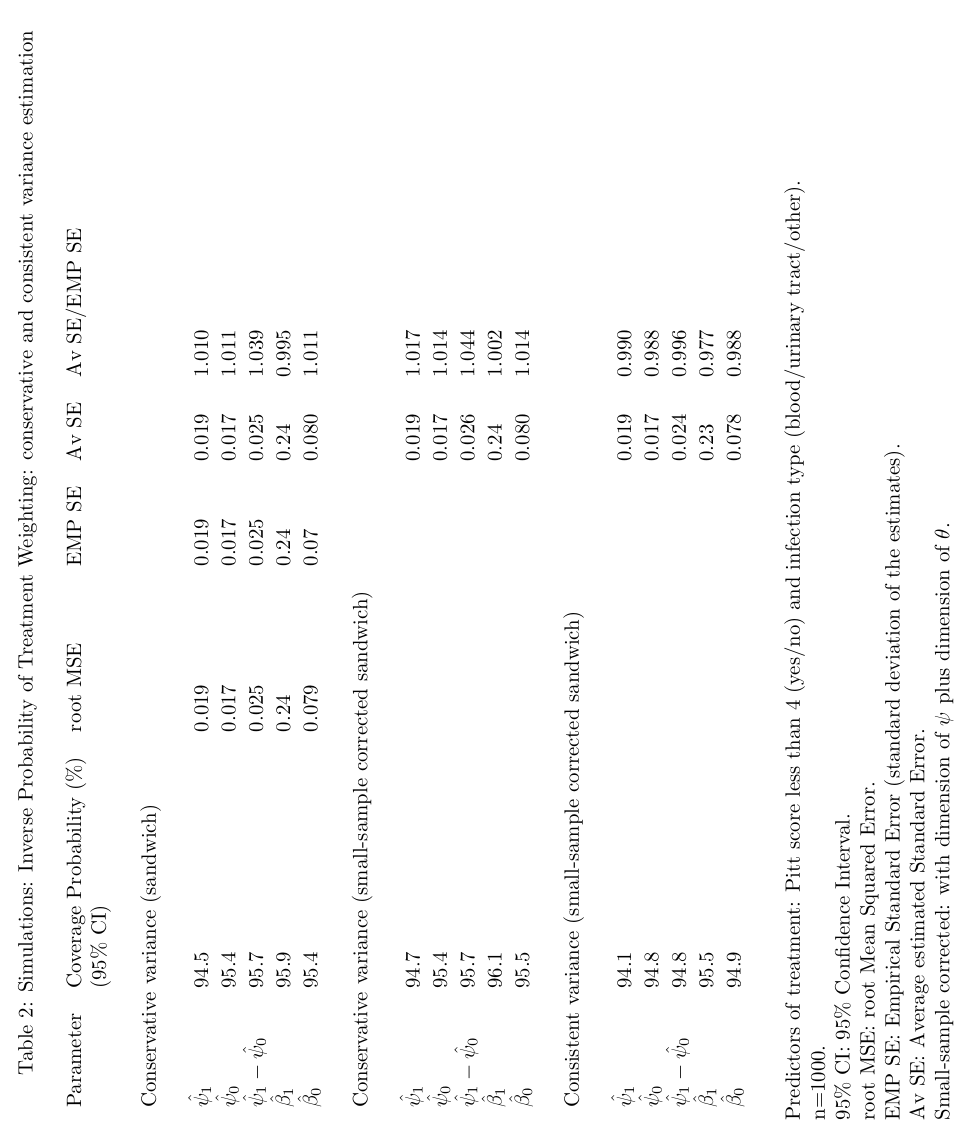}
   \includegraphics[scale=1.0,angle=0]{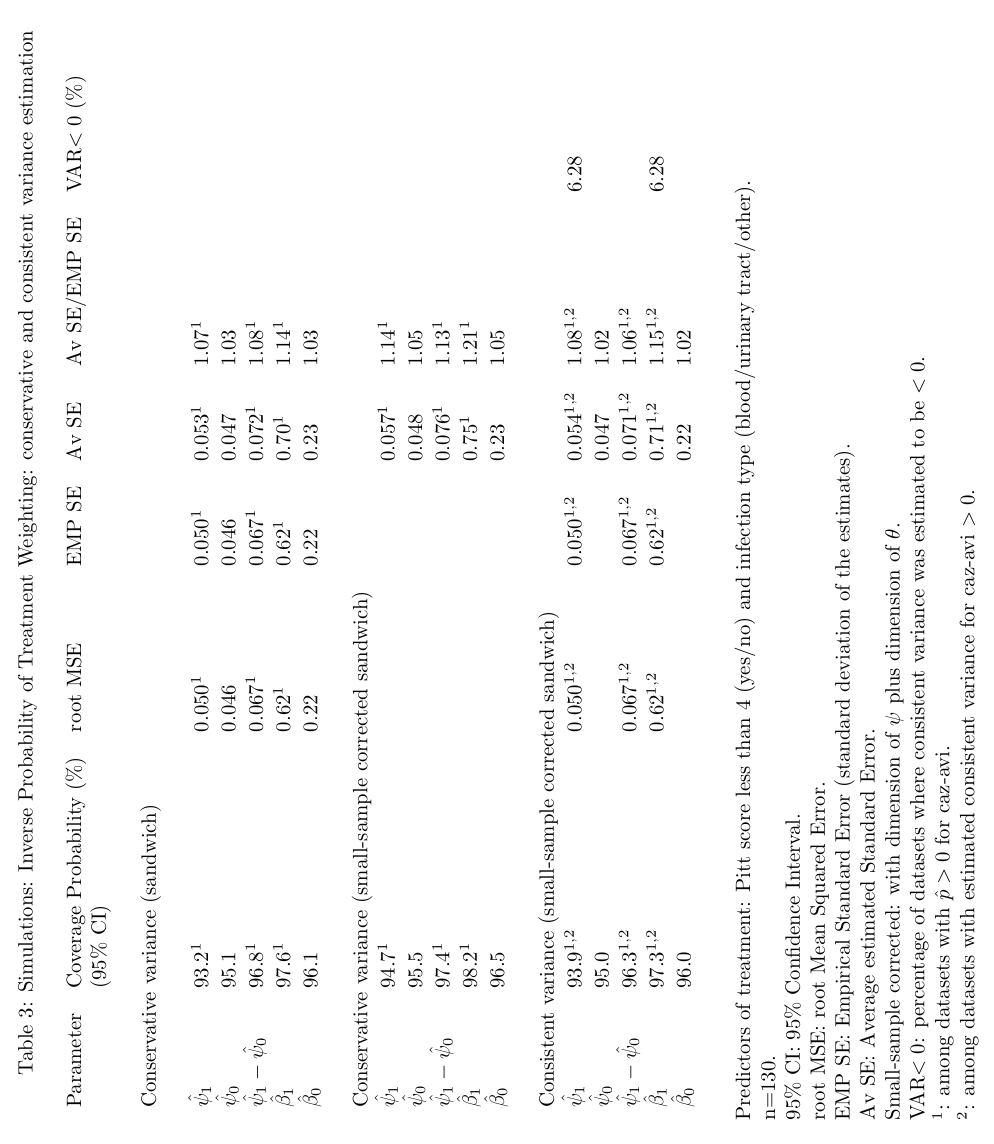}
\end{center}

If one bases the confidence intervals for the probability of hospital death on the estimates and confidence intervals for the log odds (as often done in practice) and avoids datasets that have no hospital deaths among patients treated with caz-avi (also often done in practice), all 95\% confidence intervals have reliable coverage, even for $n=130$. As expected, the consistent variance estimates have the smallest coverage probability, but very close to 95\%, and for $n=130$, possibly due to the small sample correction, over 95\%.

\section{How does the effect of ART depend on its initiation time in HIV-infected patients? Coarse Structural Nested Mean Models}\label{ART}

We estimate how the effect of one year of ART depends on the time between the estimated date of HIV infection and the ART initiation time in patients with acute and early HIV infection. This is a relevant question in HIV research, where it is known that ART improves clinical outcomes during chronic HIV infection, but less is known about the effect of ART in acute and early HIV infection. Most HIV cohorts do not include many patients in acute and early infection. The AIEDRP Core 1 database we use here includes only patients in acute and early HIV infection. 

Since interest lies in the effect of ART initiation, the time-dependent coarse Structural Nested Mean Models of Section~\ref{coarse} provide an appropriate analysis method. In the notation of Section~\ref{coarse}, we assume
\begin{equation*}
    \gamma^{(m)}_{k,\psi}\left(\overline{l}_m\right)=
\left(\psi_1+\psi_2m+\psi_3m^2\right)(k-m)1_{\left\{k>m\right\}}
\end{equation*}
as in equation~(\ref{modelex}). 
Consequently, the effect of one year of ART initiated at month $m$ is $\gamma^{(m)}_{m+12}$. \cite{Victor} and \cite{optLok} provided estimators and optimal estimators, respectively. Their confidence intervals relied on the bootstrap, and the bootstrap used to run for about 8 hours on a large computer cluster.

In this HIV application, there is both confounding by indication and censoring. To account for censoring, Inverse Probability of Censoring Weighting (\cite{RRZ}) is applied, which assumes that the missingness is Missing At Random (\cite{RubinMAR}): missingness can depend on a unit/patient's past observed information, but not on current, future, or unobserved information relevant to a patient's prognosis. With $C_{ki}=1$ if patient $i$ is censored by time $k$ and $C_{ki}=0$ if not, Missing At Random is formalized as 
\begin{assumption} Missing At Random (MAR).\label{MAR} \begin{equation*}
    C_{ki}\cip (Y_i,Z_i)|\overline{L}_{k-1i},\overline{A}_{k-1i}.
\end{equation*}
\end{assumption}
We assume that the censoring probabilities are correctly specified with a pooled logistic regression model with parameter $\theta^{(C)}$.

The nuisance parameters $\theta=(\theta^{(A)},\theta^{(C)})$ parameterize the treatment ($\theta^{(A)}$) and censoring ($\theta^{(C)}$) probabilities. As in Section~\ref{long}, one specifies a model $p_{\theta^{(A)}}(A_{k}|\overline{L}_{k},\overline{A}_{k-1})$ for the probability of treatment at time $k$ given the patient's past treatment- and covariate history. In this HIV application, we assume that $p_{\theta^{(A)}}(A_{k}|\overline{L}_{k},\overline{A}_{k-1})$ follows a pooled logistic regression model, pooled over the different time points $k$. 

Estimating $\theta^{(A)}$ with standard software for logistic regression solves the (partial) score equations (\ref{scoreMSM}) for $\theta^{(A)}$, and the estimator for the censoring parameter $\theta^{(C)}$ solves similar (partial) score equations. We chose $\vec{q}$ as in Theorem~\ref{see} approximately optimal among $\vec{q}$ with $q_k^m=0$ for $k\neq m+12$ (Section~7.2 of \cite{optLok}, \cite{Victor}). Estimating this $\vec{q}$ does not affect the asymptotic variance of $\hat{\psi}$ because the estimating equations are unbiased for any $\vec{q}$, so Section~\ref{identnotneeded} applies. Thus, because of Section~\ref{conservative}, if the pooled logistic regression models for treatment initiation and censoring are correctly specified, the sandwich estimator of the variance of $\hat{\psi}$ ignoring estimation of $\theta=(\theta^{(A)},\theta^{(C)})$ is conservative. 

Since interest lies in the effect of one year of ART, we can restrict the model assumptions on $\gamma$ to $k\leq m+12$ to avoid unnecessary model misspecification. Then, from Theorem~\ref{see} and using Inverse Probability of Censoring Weighting to account for dropout, the estimator $\hat{\psi}$ solves the linear equation
\begin{eqnarray*}
\lefteqn{P_n U_1=P_n\sum_{m=0}^{K}\sum_{k=(m+1)\vee 12}^{(m+12)\wedge (K+1)} \vec{\hat{q}}^{\;opt,k}_{m}\left(\overline{L}_{m}\right) \left(Y_k-1_{T<k} (\psi_1+\psi_2T+\psi_3T^2)(k-T)\right)\cdot }\\
&&\cdot 1_{\overline{A}_{m-1}=\overline{0}}
\left\{A_m-p_\theta\left(A_{m}=1|\overline{L}_{m}\overline{A}_{m-1}\right)\right\}\cdot\\
&&\cdot\frac{1_{C_k=0}}{\prod_{\tilde{m}=m+1}^{k} P\left(C_{\tilde{m}}=0|C_{\tilde{m}-1}=0,\overline{L}_{\tilde{m}-1},\overline{A}_{\tilde{m}-1}\right)}\\
&=&C_1\psi-C_2=0
\end{eqnarray*}
with
\begin{eqnarray*}
\lefteqn{C_1=P_n\sum_{m=0}^{K}\sum_{k=(m+1)\vee 12}^{(m+12)\wedge (K+1)} \vec{\hat{q}}^{\;opt,k}_{m}\left(\overline{L}_{m}\right)1_{T<k}(k-T)\left(\begin{array}{ccc}1 & T & T^2\end{array}\right)\cdot}\\
&&\cdot 1_{\overline{A}_{m-1}=\overline{0}}
\left\{A_m-p_\theta\left(A_{m}=1|\overline{L}_{m}\overline{A}_{m-1}\right)\right\}\cdot\\
&&\cdot\frac{1_{C_k=0}}{\prod_{\tilde{m}=m+1}^{k} P\left(C_{\tilde{m}}=0|C_{\tilde{m}-1}=0,\overline{L}_{\tilde{m}-1},\overline{A}_{\tilde{m}-1}\right)},
\end{eqnarray*}
\begin{eqnarray*}
\lefteqn{C_2=P_n\sum_{m=0}^{K}\sum_{k=(m+1)\vee 12}^{(m+12)\wedge (K+1)} \vec{\hat{q}}^{\;opt,k}_{m}\left(\overline{L}_{m}\right)Y_k\cdot}\\
&&\cdot 1_{\overline{A}_{m-1}=\overline{0}}
\left\{A_m-p_\theta\left(A_{m}=1|\overline{L}_{m}\overline{A}_{m-1}\right)\right\}\cdot\\
&&\cdot\frac{1_{C_k=0}}{\prod_{\tilde{m}=m+1}^{k} P\left(C_{\tilde{m}}=0|C_{\tilde{m}-1}=0,\overline{L}_{\tilde{m}-1},\overline{A}_{\tilde{m}-1}\right)}.
\end{eqnarray*}
Thus, the conservative sandwich estimator of the variance of $(\hat{\psi}_1,\hat{\psi}_2,\hat{\psi}_3)$ ignoring estimation of $\theta$ as in Equation~(\ref{sandwich}) equals
\begin{equation}\label{sandwichopt}\frac{1}{n}C_1^{-1}DC_1^{-1\top} \;\;\;\;\;\;\;\;\;\;\;\;\;\;\text{ with }
D= P_n U_1^{\otimes 2}(\hat{\psi},\hat{\theta}).
\end{equation}

Because of Assumption~\ref{MAR}, the estimating equations for the censoring model $U_{2}^{(C)}$ are orthogonal to the estimating equations for the treatment initiation model $U_{2}^{(A)}$; this follows from the Law of Iterated Expectations, conditioning on $(Y,Z)$, which fixes everything except for $C_k$, which has conditional expectation $P\left(C_k=0|C_{k-1}=0, \overline{L}_{k-1}, \overline{A}_{k-1}\right)$ (details in Web-Appendix~\ref{CAorthogonal}). Thus, 
\begin{equation*}
    E\left(U_2^{\otimes 2}\right)=E\left(\left(\begin{array}{c} U_{2}^{(A)}\\U_{2}^{(C)}\end{array}\right)^{\otimes 2}\right)
    =\left(\begin{array}{cc} E\left(U_{2}^{(A)\otimes 2}\right) & 0 \\0 & E\left(U_{2}^{(C)\otimes 2}\right)\end{array}\right).
\end{equation*}
Then, from Section~\ref{varcorrection} equation~(\ref{varcorterm}), the variance can be consistently estimated by the sandwich estimator ignoring estimation of $\theta$, (\ref{sandwichopt}), minus an estimate of
\begin{eqnarray*} \lefteqn{C_1^{-1}\left(EU_1U_{2}^{(A)\top} \;  EU_1U_{2}^{(C)\top}\right)\left(\begin{array}{cc} \widehat{VAR}\,\theta^{(A)} & 0\\ 0 & \widehat{VAR}\,\theta^{(C)}\end{array}\right) \left(EU_1U_{2}^{(A)\top} \; EU_1U_{2}^{(C)\top}\right)^\top C_1^{-1\,\top}}\\
&=&C_1^{-1}\bigl(EU_1U_{2}^{(A)\top}\bigr) \widehat{VAR}\,\theta^{(A)} \bigl(EU_1U_{2}^{(A)\top}\bigr)^\top C_1^{-1\top}\\
&&+C_1^{-1}\bigl(EU_1U_{2}^{(C)\top}\bigr) \widehat{VAR}\,\theta^{(C)} \bigl(EU_1U_{2}^{(C)\top}\bigr)^\top C_1^{-1\top},\end{eqnarray*}
where each expectation can be estimated by its empirical average evaluated at the estimates $\hat{\psi}$ and $\hat{\theta}=(\hat{\theta}^{(A)},\hat{\theta}^{(C)})$.

\begin{center}
   \includegraphics[scale=0.9,angle=0]{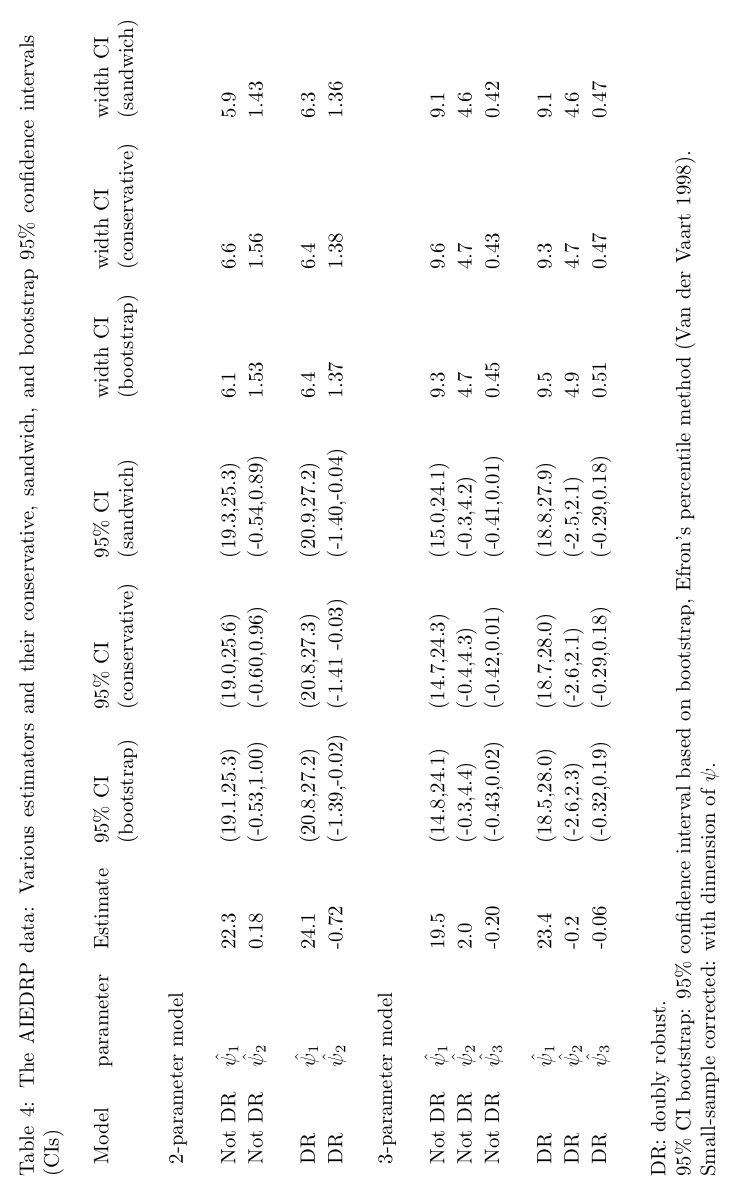}
\end{center}

Table~4 describes results for 2 estimators for the 3-dimensional and 2-dimensional (omitting $\psi_3$) parameter $\psi$ from Equation~(\ref{modelex}). The first estimator is not doubly robust, the second estimator is doubly robust as outlined in \cite{optLok}: for consistency it depends on an outcome regression model \emph{or} the treatment initiation model, but it does require correct specification of the censoring model. Table~2 reports the (conservative) sandwich estimates of the variance ignoring estimation of $\theta$, the (consistent) sandwich estimates of the variance using the correction terms, and the (consistent) bootstrap estimates of the variance. We see that as expected, the bootstrap and the consistent variance estimate are reasonably close. Perhaps somewhat surprisingly, in this application the conservative confidence intervals are in many cases closer to the bootstrap confidence intervals than those based on the consistent sandwich estimator of the variance.

Figure~1 compares bootstrap confidence intervals (Efron's percentile method, see e.g.~\cite{Vaart}) with conservative (Figure~1a) and consistent (Figure~1b) sandwich estimator-based confidence intervals (based on our theory) for the effect of 1 year of ART initiated at month $m$ since infection (with month on the x-axis). These sandwich estimators  use that
\begin{equation*}
    \widehat{VAR}(\hat{\psi}_1+\hat{\psi}_2m)=\widehat{VAR}(\hat{\psi}_1)+2m\widehat{COV}(\hat{\psi}_1,\hat{\psi}_2)+m^2\widehat{VAR}(\hat{\psi}_2).
\end{equation*}

\begin{center}
   \includegraphics[scale=1,angle=0]{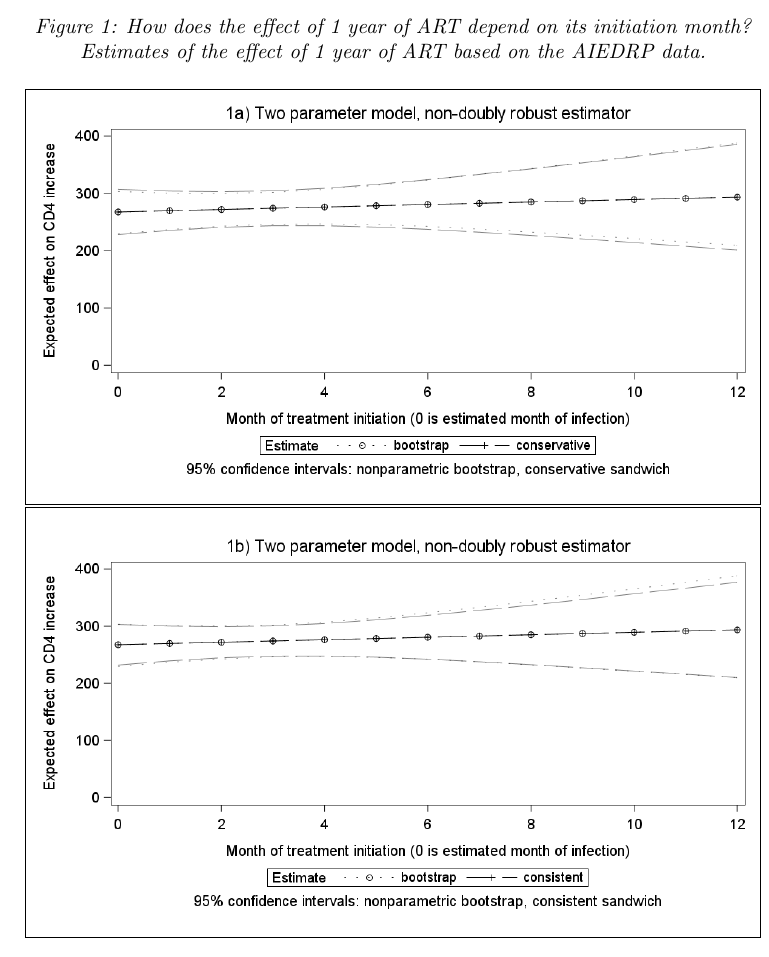}
\end{center}


\section{Discussion}\label{Discussion}

This paper has 5 main results for estimators $\hat{\psi}$ based on unbiased estimating equations including a nuisance parameter $\theta$. First, it provides a general formula for the variance when nuisance parameters are estimated by solving unbiased estimating equations. In general, using the Partition Inverse Formula, it turns out to be fairly easy to avoid the bootstrap.
That being said, Section~\ref{bootstrap} proves that the bootstrap, while for large datasets often requiring considerable computing time, often leads to consistent variance estimates as well.

The other 4 results apply when the nuisance parameter $\theta$ is estimated by solving (partial) score equations and the parameter of interest $\psi$ does not depend on $\theta$: 1.\ The true limiting variance of $\hat{\psi}$ is smaller or remains the same when the nuisance parameter $\theta$ is estimated, compared to if $\theta$ were known and plugged in. 2.\ If estimation of the nuisance parameter $\theta$ is ignored, the resulting sandwich estimator for the variance of $\hat{\psi}$ is conservative. 3.\ We provide a variance correction. 4.\ We show that if the estimator $\hat{\psi}$ for the parameter of interest $\psi$ with the true nuisance parameter $\theta$ plugged in is efficient, the limiting variance of $\hat{\psi}$ does not depend on whether or not $\theta$ is estimated, and the sandwich estimator for the variance of $\hat{\psi}$ ignoring estimation of $\theta$ is consistent.

The data applications in Sections~\ref{cazavi} and~\ref{ART} and the simulations in Section~\ref{simulations} indicate that the conservative estimator of the variance, ignoring estimation of the nuisance parameter, works well in these practical applications. The simulations in Section~\ref{simulations} indicate that the variance estimator correcting for estimation of the nuisance parameters work well in larger samples. For small samples, even with a small sample correction, it seems better to use the bootstrap or the conservative estimator of the variance ignoring estimation of the nuisance parameter. In any case, it is more important to avoid the bootstrap for larger samples, since its main disadvantage is computing time, which is not as much of an issue in smaller samples. 

The examples in this article show that the sandwich estimator ignoring estimation of the nuisance parameter seems to often work well when the nuisance parameter is estimated by solving partial score equations and the parameter of interest does not depend on the nuisance parameter. The examples showed limited gains from correcting this conservative estimator of the variance. One can even expect some upward bias of the correction term, which by its form is always positive semi-definite. More examples and simulations are needed to see when one can expect larger gains from correction of the conservative variance.

The sandwich estimator for the variance of $\hat{\psi}$ does not need to be conservative if the nuisance parameter model is misspecified. In such settings the parameter of interest $\psi$ may not even be identified. Doubly robust methods, see Section~\ref{DR}, can be used to reduce the impact of model misspecification.

The sandwich estimator for the variance of $\hat{\psi}$ also does not need to be conservative if Assumption~\ref{psinotftheta} does not hold because the parameter of interest $\psi$ depends on the nuisance parameter $\theta$. An example is the estimator of the mean under say treatment $1$ from observational data based on a conditioning argument (see e.g.~\cite{rubin1977assignment}). There, in the notation of Section~\ref{cazavi}, one first estimates $E_\theta[Y|A=1,L]$, and then
\begin{equation*}
    \hat{E}Y^{(1)}=\frac{1}{n}\sum_{i=1}^n\hat{E}_{\hat{\theta}}[Y|A=1,L=l_i].
\end{equation*}
If $Y$ follows a regression model $Y_i=\theta_0+\theta_1L_i+\epsilon_i$ with $E[\epsilon_i|L_i]=0$, a conditioning argument shows that $U_1$ and $U_2$ are orthogonal. Then, from e.g.\ equation~(\ref{VARhatgenest}), it follows that inference on $\psi$ ignoring estimation of $\theta$ is anti-conservative. In this setting, $\psi=\theta_0+\theta_1E(L)$, and equation~(\ref{derivative0}) does not hold, since its derivation used that the derivative with respect to $\theta$ did not involves $\psi'(\theta)$. Equation~(\ref{VARhatgenest}) can still be used to estimate the variance consistently in such settings.

\cite{henmi2004paradox} developed similar theory for settings where the likelihood is parameterized by the parameter of interest $\psi$, the finite-dimensional nuisance parameter $\theta$, and an infinite-dimensional nuisance parameter $\eta$, which all vary independently. Our illustrative examples do not follow this structure. For Inverse Probability of Treatment Weighting (Sections~\ref{IPTWpoint}, , \ref{IPTWtimedep}, and~\ref{cazavi}), $\theta$ parameterizes what we see given baseline (or time-varying) covariates, and any likelihood factorization involving $\theta$ will include outcome distributions given those covariates, whereas interest often lies in unconditional means $\psi$ (or in the case of Marginal Structural Models sometimes means given baseline covariates); then $\psi$ does not appear directly in the factorized likelihood. In the case of time-dependent coarse Structural Nested Models (Sections~\ref{coarse} and~\ref{ART}), $\theta$ parameterizes what we see given time-varying covariates, but interest lies in a parameter $\psi$ describing means conditional on covariates measured (much) earlier than the time the outcome is measured; also here, $\psi$ does not appear directly in the factorized likelihood. There are however interesting settings where the conditions of \cite{henmi2004paradox} are satisfied. Their first illustrative example describes a setting from \cite{robins1992estimating} where (in the notation of Section~\ref{IPTWpoint}) $Y_i^{(a)}=\psi a + h(L_i)+\epsilon^{(a)}_i$ with $E[\epsilon^{(a)}_i|A_i=a,L_i]=0$, and where $A_i$ follows a logistic regression model. This constant treatment effect model (given $L$) renders a natural parameterization of the likelihood factorizing into a factor involving $\psi$ and a factor involving $\theta$. In their Inverse Probability of Censoring Weighting example, \cite{henmi2004paradox} assume that $(Y,L)$ follow a joint normal distribution with $EY=\psi$, $EL=0$, and VAR$L=1$. This leads to a parameterization of $Y$ given $L$ parameterized by $\psi$ and a factorized fully parametric likelihood. In contrast, in many Inverse Probability of Censoring and Inverse Probability of Treatment Weighting analyses, no assumptions are made on how $Y$ and $L$ are related. In these and many other settings,  $\psi$ is the solution to an identifying equation that depends on the unknown data distribution, and the likelihood does not explicitly depend on $\psi$; Sections~\ref{Examples}, \ref{cazavi}, and \ref{ART} provide examples.

Interesting areas for future research include variance estimation after model selection for the nuisance parameter model. In addition, Lok, Buchanan, and Spiegelman worked out the details of variance derivations for Inverse Probability of Treatment Weighting of Marginal Structural Models with time-dependent treatment at the cluster level, with observations not necessarily independent. This work will soon be posted on arxiv. Another interesting area for future research is variance estimation when the nuisance parameter is not estimated using unbiased estimating equations.

The examples in this paper come from causal inference, but the results included in this paper are much more general: the main result (\ref{VARhatgen}) applies to any setting where nuisance parameters are estimated using unbiased estimating equations, and the main result (\ref{g-var}) applies to any setting where the parameter of interest does not depend on the nuisance parameter, and the nuisance parameter is estimated by solving (partial) score equations. In practice, for larger datasets, one could first calculate the confidence intervals based on the conservative estimator of the variance, and if the confidence intervals are not satisfactory, one can either use the bootstrap and wait for the results, or base confidence intervals on the consistent estimator of the variance.


\section*{Acknowledgements} I am extremely grateful to the participants with bacterial infections in the ARLG Crackle~1 study, and to the HIV-infected participants who volunteered for the AIEDRP study. I thank David van~Duin, Scott Evans, and Lauren Komarov for granting access to the ARLG Crackle~1 data used in Section~\ref{cazavi}, and for a fruitful collaboration on \cite{Duin}. I thank Victor DeGruttola, Davey Smith, and Susan Little for granting access to the AIEDRP Core~01 data used in Section~\ref{ART}, and for a fruitful collaboration on \cite{Victor}. I thank Andrea Rotnitzky and Oliver Dukes for helpful comments and suggestions. 


This research was supported by NSF DMS 1854934 to Lok, ask about ARLG, ask about AIEDRP. The content is solely the responsibility of the author and does not necessarily represent the official views of the National Institutes of Health or the National Science Foundation.

\addcontentsline{toc}{chapter}{Bibliography}
\bibliographystyle{chicago} \bibliography{ref}

\appendix

\section{Proofs}
\label{Proofs}

\noindent {\bf Proof of Equation~(\ref{IPTWwithtime}).}
\begin{eqnarray*}
\lefteqn{E\left(U_1(\psi_1^*,\theta^*)\right)}\\
&=&\sum_{\overline{a}_K}E\left( \frac{1_{\overline{A}_{Ki}=\overline{a}_{K}}}{\prod_{k=0}^K P\left(A_{ki}=a_{ki}|\overline{L}_{ki},\overline{A}_{k-1i}=\overline{a}_{k-1}\right)}\left(Y_i^{(\overline{a}_K)}-\mu_{\psi^*}^{(\overline{a}_K)}\right)
\right)\\
&=&\sum_{\overline{a}_K}E\left(E\left[\frac{1_{\overline{A}_{Ki}=\overline{a}_{K}}}{\prod_{k=0}^K P\left(A_{ki}=a_{k}|\overline{L}_{ki},\overline{A}_{k-1i}=\overline{a}_{k-1}\right)}\left(Y_i^{(\overline{a}_K)}-\mu_{\psi^*}^{(\overline{a}_K)}\right)\mid Y_i^{(\overline{a}_K)},\overline{L}_{Ki}^{(\overline{a}_K)}\right]\right)\\
&=&\sum_{\overline{a}_K}E\left(\frac{E\left[1_{\overline{A}_{Ki}=\overline{a}_{K}}\mid Y_i^{(\overline{a}_K)},\overline{L}_{Ki}^{(\overline{a}_K)}\right]}{\prod_{k=0}^K P\left(A_{ki}=a_{k}|\overline{L}_{ki},\overline{A}_{k-1i}=\overline{a}_{k-1}\right)}\left(Y_i^{(\overline{a}_K)}-\mu_{\psi^*}^{(\overline{a}_K)}\right)\right)\\
&=&\sum_{\overline{a}_K}E\left(\frac{E\left[1_{\overline{A}_{Ki}=\overline{a}_{K}}\mid \overline{L}_{Ki}^{(\overline{a}_K)}\right]}{\prod_{k=0}^K P\left(A_{ki}=a_{ki}|\overline{L}_{ki},\overline{A}_{k-1i}=\overline{a}_{k-1}\right)}\left(Y_i^{(\overline{a}_K)}-\mu_{\psi^*}^{(\overline{a}_K)}\right)\right)\\
&=&\sum_{\overline{a}_K}E\left(Y_i^{(\overline{a}_K)}-\mu_{\psi^*}^{(\overline{a}_K)}\right)=0,
\end{eqnarray*}
where for the first equality we use that on $\overline{A}_K=\overline{a}_K$, we observe the outcome $Y^{(\overline{a}_K)}$, for the second equality we use the Law of Iterated Expectations, for the fourth equality we use the Assumption of No Unmeasured Confounding~\ref{nucwithtime}, and in the last line we used that on $\overline{A}_K=\overline{a}_K$,
\begin{eqnarray*}
P\left(\overline{A}_{Ki}=\overline{a}_{K}\mid \overline{L}_{Ki}^{(\overline{a}_K)}\right)
&=&\prod_{k=0}^K P\left(A_{ki}=a_{k}|\overline{L}_{Ki}^{(\overline{a}_K)},\overline{A}_{k-1i}=\overline{a}_{k-1}\right)\\
&=&\prod_{k=0}^K P\left(A_{ki}=a_{k}|\overline{L}_{ki},\overline{A}_{k-1i}=\overline{a}_{k-1}\right)
\end{eqnarray*}
because of the Assumption of No Unmeasured Confounding~\ref{nucwithtime} and since on $\overline{A}_K=\overline{a}_K$, $\overline{L}_{ki}=\overline{L}_{ki}^{(\overline{a}_K)}$.

It is easy to see that the same proof works if we include the stabilizer of the weights in the numerator, except that if $X_i$ is included in the stabilizer, the stabilizer of the weights is still there in the last line, and we need to condition on $X_i$ to show that this last expectation is $0$.

\begin{lemma}\label{Donskerconv} If $\hat{\psi}$ is a consistent estimator of a parameter $\psi$, $X_i$ for $i=1,\ldots,n$ are independent and identically distributed, and $U(X_i,\psi)$ is a function of $\psi$ that is bounded by an integrable function and continuously differentiable with a derivative that is bounded by a square integrable function, then
\begin{equation*}
P_n(U(X_i,\hat{\psi})) \rightarrow^P E\left(U(X_i,\psi^*)\right).
\end{equation*}
\end{lemma}

\noindent {\bf Proof:} $U(X_i,\psi)$ form a Donsker class by Example~19.7 in \cite{Vaart}, so that
\begin{equation*}
\sup_{\psi}\left\|P_n U(X_i,\psi)- E\left(U(X_i,\psi)\right)\right\|\rightarrow^P 0.
\end{equation*}
Moreover, since $U(X_i,\psi)$ is continuously differentiable in $\psi$, $U(X_i,\hat{\psi})\rightarrow^P U(X_i,\psi^*)$, and since it is also bounded by an integrable function, Lebesgue's Dominated Convergence Theorem implies that 
$E\big(U(X_i,\hat{\psi})\bigr)$
converges to $E\left(U(X_i,\psi^*)\right)$. Combining, and using the triangle inequality, we find that indeed Lemma~\ref{Donskerconv} holds.
\hfill $\Box$\bigskip

\begin{lemma}\label{pFisher}
Under Assumptions~\ref{diffb} and~\ref{support},  $E\frac{\partial}{\partial \theta}U_2=-EU_2^{\otimes 2}$. 
\end{lemma}

\noindent {\bf Proof} This well-known fact is proven for the partial score the same way it is proven for the score:
\begin{eqnarray*}
E\frac{\partial}{\partial\theta} U_2
&=&\int\frac{\partial}{\partial \theta}\left(\frac{\partial}{\partial \theta} \log f_\theta\right) f_\theta dz\\
&=&\frac{\partial}{\partial \theta}\int\left(\frac{\partial}{\partial \theta} \log f_\theta\right) f_\theta dz-\int\left(\frac{\partial}{\partial \theta} \log f_\theta\right) \frac{\partial}{\partial \theta}f_\theta dz\\
&=&\frac{\partial}{\partial \theta}EU_2-\int\left(\frac{\partial}{\partial \theta} \log f_\theta\right)^{\otimes 2} f_\theta dz\\
&=&0-EU_2^{\otimes 2}=-EU_2^{\otimes 2},
\end{eqnarray*}
where for the second equality we used Assumption~\ref{support}.\\

\noindent {\bf Proof} of Equation~(\ref{ddtheta}).
\begin{eqnarray*}
\lefteqn{E\left.\frac{\partial}{\partial \theta}\right|_{\theta^*}U_1(\psi^*,\theta)}\\
&=&E\left.\frac{\partial}{\partial \theta}\right|_{\theta^*} \tilde{U}_1(\psi^*,\theta)+A_{\psi^*,\theta^*}E\left.\frac{\partial}{\partial \theta}\right|_{\theta^*} U_2(\theta)\\
&=&A_{\psi^*,\theta^*}E\left.\frac{\partial}{\partial \theta}\right|_{\theta^*} U_2(\theta)\\
&=&-E\left(U_1(\psi^*,\theta^*)U_2(\theta^*)^\top\right)\left(E\left(U^{\otimes 2}_2(\theta^*)\right)\right)^{-1}E\left(U^{\otimes 2}_2(\theta^*)\right)\\
&=&-EU_1U_2^\top,
\end{eqnarray*}
where for the first equality we use equation~(\ref{U1tilde}), for the second equality we use equation~(\ref{derivative0}) and for the third equality we use equation~(\ref{Apsitheta}) and that $U_2$ is a partial score, so that under standard regularity conditions (e.g.\ Assumptions~\ref{diffb} and~\ref{support}, see Lemma~\ref{pFisher}), $E\frac{\partial}{\partial \theta}U_2=-EU_2^{\otimes 2}$. 

Equation~(\ref{ddtheta}) was shown for the specific case of Inverse Probability of Censoring Weighting in \cite{RR95} (proof of their Proposition~5, referencing also \cite{pierce1982asymptotic}).
\hfill $\Box$\bigskip

\noindent {\bf Proof} of expression (\ref{g-asvar}) for the limiting variance. 
\begin{eqnarray*}
\lefteqn{E\left(\left(-U_{1}+E(U_{1}U_{2}^\top)
E\left(U_{2}^{\otimes 2}\right)^{-1}U_{2}\right)\left(-U_{1}+E(U_{1}U_{2}^\top)
E\left(U_{2}^{\otimes 2}\right)^{-1}U_{2}\right)^\top\right)}\\
&=&E\left(U_{1}^{\otimes 2}\right)
-E\left(U_{1}U_{2}^\top\right)
E\left(U_{2}^{\otimes 2}\right)^{-1}\left(E\left(U_{1}U_{2}^\top\right)\right)^\top\\
&&-E\left(U_{1}U_{2}^\top\right)E\left(U_{2}^{\otimes 2}\right)^{-1}(E(U_{1}U_{2}^\top))^\top\\
&&+E\left(U_{1}U_{2}^\top\right) E\left(U_{2}^{\otimes 2}\right)^{-1}E\left(U_{2}^{\otimes 2}\right)E\left(U_{2}^{\otimes 2}\right)^{-1}\left(E\left(U_{1}U_{2}^\top\right)\right)^\top\\
&=&E\left(U_{1}^{\otimes 2}\right)
-2E\left(U_{1}U_{2}^\top\right)
E\left(U_{2}^{\otimes 2}\right)^{-1}\left(E\left(U_{1}U_{2}^\top\right)\right)^\top\\
&&+E\left(U_{1}U_{2}^\top\right)E\left(U_{2}^{\otimes 2}\right)^{-1}\left(E\left(U_{1}U_{2}^\top\right)\right)^\top\\
&=&E\left(U_{1}^{\otimes 2}\right)
-E\left(U_{1}U_{2}^\top\right)
E\left(U_{2}^{\otimes 2}\right)^{-1}\left(E\left(U_{1}U_{2}^\top\right)\right)^\top.
\end{eqnarray*}

\section{Estimating the nuisance parameter $\theta$ and the variance of efficient estimators of $\psi$}\label{Efficient}

Efficient estimating equations for $\psi$ when $\theta$ is known are orthogonal to estimating equations for $\theta$. This known fact can be understood as follows. When $\theta$ is known, we can estimate $\psi$ by $\tilde{\psi}$ solving $P_n U_1^{eff}(\tilde{\psi},\theta^*)=0$. Since we assumed that $U_1^{eff}$ leads to an efficient estimator, $\tilde{\psi}$ must be efficient. The projection of $U^{eff}_1$ onto $U_2$ can be calculated as in (\ref{Apsitheta}). We estimate $C^*=E(U_2(\theta^*)^{\otimes 2})$ by $\hat{C}=P_nU_2(\theta^*)^{\otimes 2}$. We estimate 
$B^*=E(U^{eff}_1(\psi^*,\theta^*)U_2(\theta^*))$ by $\hat{B}=P_n U^{eff}_1(\tilde{\psi},\theta^*)U_2(\theta^*)$. Then, we estimate the estimating equations, the projection of $U_1^{eff}$ onto the orthocomplement of the space spanned by $U_2$, and solve those to obtain $\hat{\psi}$.

This procedure is equivalent to solving the stacked estimating equations
\begin{equation*}
P_n\left(\begin{array}{c}U^{eff}_1(\psi,\theta^*)- B C^{-1} U_2(\theta^*)\\
U^{eff}_1(\tilde{\psi},\theta^*)U_2^\top(\theta^*)-B\\
U^{eff}_1(\tilde{\psi},\theta^*)\\
U_2(\theta^*)^{\otimes 2}-C
\end{array}\right)=\vec{0}.
\end{equation*}
Because of \cite{Vaart} Theorem~5.21, under regularity conditions
\begin{equation}\label{thetastar}
\left(\begin{array}{c}
\hat{\psi}-\psi^*\\
\hat{B}-B^*\\
\hat{\tilde{\psi}}-\psi^*\\
\hat{C}-C^*
\end{array}\right)
=D(\psi^*,\theta^*)^{-1}
P_n\left(\begin{array}{c}U^{eff}_1(\psi^*,\theta^*)- B C^{-1} U_2(\theta^*)\\
U^{eff}_1(\psi^*,\theta^*)U_2^\top(\theta^*)-B\\
U^{eff}_1(\psi^*,\theta^*)\\
U_2(\theta^*)^{\otimes 2}-C
\end{array}\right)+o_P(1/\sqrt{n}),
\end{equation}
with 
\begin{equation*}
D(\psi^*,\theta^*)=\left(\begin{array}{cccc}
E\frac{\partial}{\partial \psi} U^{eff}_1 & 0 & 0 & 0\\
0 & -I & E \frac{\partial}{\partial \psi} U^{eff}_1 U_2^\top & 0\\
0 & 0 & E\frac{\partial}{\partial \psi} U^{eff}_1 & 0\\
0 & 0 & 0 & -I
\end{array}\right).
\end{equation*}
To show that the first row of
\begin{equation*}
D(\psi^*,\theta^*)=E\left(\frac{\partial}{\partial (\psi,B,\tilde{\psi},C)}\left(\begin{array}{c}U^{eff}_1(\psi,\theta^*)- B C^{-1} U_2(\theta^*)\\
U^{eff}_1(\tilde{\psi},\theta^*)U_2^\top(\theta^*)-B\\
U^{eff}_1(\tilde{\psi},\theta^*)\\
U_2(\theta^*)^{\otimes 2}-C
\end{array}\right)\right)
\end{equation*}
has this form can be shown as follows. The derivative of the first entry of the estimating equation with respect to $B$ and $C$ has expectation $0$ since it has the factor $U_2(\theta^*)$, which has expectation $0$, and the first entry does not involve $\tilde{\psi}$, so the derivative with respect to $\tilde{\psi}$ is equal to $0$.

$D(\psi^*,\theta^*)$ has an inverse under Regularity Condition~\ref{inverses}. It is easy to see that the first row of $D(\psi^*,\theta^*)^{-1}$ needs to have  $(E\frac{\partial}{\partial \psi} U^{eff}_1)^{-1}$ in the first column and zeros in the next columns. Combining this with (\ref{thetastar}) implies that
\begin{equation}
\hat{\psi}-\psi^*=(E\frac{\partial}{\partial \psi} U^{eff}_1)^{-1} P_n\left(U^{eff}_1(\psi^*,\theta^*)- B C^{-1} U_2(\theta^*)\right) + o_P(1/\sqrt{n}).
\end{equation}
$U^{eff}_1(\psi^*,\theta^*)- B C^{-1} U_2(\theta^*)$ is the projection of $U_1^{eff}$ onto the orthocomplement of the space spanned by $U_2$, and so if the projection changes $U^{eff}_1$, $\hat{\psi}$ is more efficient than $\hat{\tilde{\psi}}$, which was efficient. Thus, the projection does \emph{not} change $U^{eff}_1$, and $U^{eff}_1$ is already orthogonal to $U_2$.\bigskip

If the estimating equations $U^{eff}_1$ for $\psi$ are efficient were $\theta$ known, and in addition 1.\ $\psi$ does not depend on $\theta$ and 2.\ $\theta$ is estimated by solving (partial) score equations, estimation of $\theta$ does not affect the variance of $\hat{\psi}$. The proof of this statement from Section~\ref{efficient} is as follows. Because of the beginning of Web-appendix~\ref{Efficient}, $U^{eff}_1$ is orthogonal to the partial score $U_2$. From Section~\ref{replace} it follows that $\tilde{U}^{eff}_1=U^{eff}_1$. Then, we know from Section~\ref{nuisestVAR} that 1.\ estimating the nuisance parameter $\theta$ does not affect the variance of the estimator of $\psi$, and 2.\ the sandwich estimator for the variance of the estimator for $\psi$ which estimates $\theta$ in a first step by solving the partial score equations is consistent.

\section{The treatment initiation and the censoring model lead to orthogonal estimating equations}\label{CAorthogonal}

Because of Assumption~\ref{MAR}, the estimating equations for the censoring model $U_{2}^{(C)}$ are orthogonal to the estimating equations for the treatment initiation model $U_{2}^{(A)}$. This can be seen as follows.
\begin{eqnarray*}\lefteqn{E\left(U_2^{(A)}U_2^{(C)\top}\right)}\\
&=&E\left(\sum_{k=0}^K \sum_{m=1}^K E\left[\left(\begin{array}{c} L_{ki}\\ A_{k-1i} \end{array}\right) (A_{ki}-p_{\theta^{(A)}}(A_{ki}=1|\overline{L}_{ki},\overline{A}_{k-1i})\right.\right.\\
&&\cdot \left.\left.\left(\begin{array}{c} L_{m-1i}\\ A_{m-1i} \end{array}\right)^\top 1_{C_{m-1i}=0}(C_{mi}-p_{\theta^{(C)}}\left(C_{mi}=1|\overline{L}_{m-1i},\overline{A}_{m-1i},C_{m-1i}=0\right))|\overline{A}_K,\overline{L}_K,C_{m-1i}\right]\right)\\
&=&E\left(\sum_{k=0}^K \sum_{m=1}^K \left(\begin{array}{c} L_{ki}\\ A_{k-1i} \end{array}\right) (A_{ki}-p_{\theta^{(A)}}(A_{ki}=1|\overline{L}_{ki},\overline{A}_{k-1i})\right.\\
&&\cdot \left.\left(\begin{array}{c} L_{m-1i}\\ A_{m-1i} \end{array}\right)^\top 1_{C_{m-1i}=0}( E\left[C_{mi}|\overline{A}_K,\overline{L}_K,C_{m-1i}=0\right]-p_{\theta^{(C)}}\left(C_{mi}=1|\overline{L}_{m-1i},\overline{A}_{m-1i},C_{m-1i}=0\right))\right)\\
&=&0.
\end{eqnarray*}
We conditioned $\overline{A}_K,\overline{L}_K,C_{m-1i}$, using the Law of Iterated Expectations, which fixes everything except for $C_m$, which has conditional expectation $P\left(C_m=0|\overline{L}_{m-1}, \overline{A}_{m-1}, C_{m-1}=0\right)$. 

\section{Time-dependent coarse Structural Nested Mean Models}

We show that Assumption~\ref{intconf2} implies that for $k>m$,
\begin{equation}\label{gammaalt}
\hspace*{3cm}\gamma_k^{(m)}\left(\overline{l}_m\right)=E\left[Y_{k}^{(m)}-Y_{k}^{(\emptyset)}\mid \overline{L}^{(\emptyset)}_{m}=\overline{l}_m\right],
\end{equation}
so that Assumption~\ref{psinotftheta} is satisfied.

To prove equation~(\ref{gammaalt}), we show that under Assumption~\ref{intconf2},
\begin{equation}\label{coarsealt}
    P\left(Y_{k}^{(m)}=y\mid \overline{L}^{(\emptyset)}_{m}=\overline{l}_m,T=m\right)=P\left(Y_{k}^{(m)}=y\mid \overline{L}^{(\emptyset)}_{m}=\overline{l}_m\right),
\end{equation}
and similarly for $Y_{k}^{(\emptyset)}$. Equation~(\ref{gammaalt}) follows. First, notice that
\begin{eqnarray*}
P\left(Y_{k}^{(m)}=y\mid \overline{L}^{(\emptyset)}_{m}=\overline{l}_m,T=m\right)
&=&
P\left(Y_{k}^{(m)}=y\mid \overline{L}^{(\emptyset)}_{m}=\overline{l}_m,\overline{A}_{m-1}=\overline{0},A_m=1\right)\\
&=&P\left(Y_{k}^{(m)}=y\mid \overline{L}^{(\emptyset)}_{m}=\overline{l}_m,\overline{A}_{m-1}=\overline{0}\right),\end{eqnarray*}
where for the second equality we used Assumption~\ref{intconf2}. Next, continue by backwards induction. We show that if \begin{equation*}
    P\left(Y_{k}^{(m)}=y\mid \overline{L}^{(\emptyset)}_{m}=\overline{l}_m,T=m\right)=P\left(Y_{k}^{(m)}=y\mid \overline{L}^{(\emptyset)}_{m}=\overline{l}_m,\overline{A}_p=\overline{0}\right)
\end{equation*}
holds for some $p\leq m-1$, it also holds for $p-1$. The result for $m-1$ was shown above. Equation~(\ref{coarsealt}) follows from the final result for $p=-1$ (no conditioning on $A$). The induction step:
\begin{eqnarray*}
P\left(Y_{k}^{(m)}=y\mid \overline{L}^{(\emptyset)}_{m}=\overline{l}_m,T=m\right)
&=&P\left(Y_{k}^{(m)}=y\mid \overline{L}^{(\emptyset)}_{m}=\overline{l}_m,\overline{A}_{p}=\overline{0}\right)\\
&=&\frac{P\left(Y_{k}^{(m)}=y,\overline{L}^{(\emptyset)}_{m}=\overline{l}_m,\overline{A}_{p}=\overline{0},\overline{L}^{(\emptyset)}_{p}=\overline{l}_{p}\right)}{P\left(\overline{L}^{(\emptyset)}_{m}=\overline{l}_m,\overline{A}_{p}=\overline{0},\overline{L}^{(\emptyset)}_{p}=\overline{l}_{p}\right)}\\
&=&\frac{P\left(Y_{k}^{(m)}=y,\overline{L}^{(\emptyset)}_{m}=\overline{l}_m\mid \overline{A}_{p}=\overline{0},\overline{L}^{(\emptyset)}_{p}=\overline{l}_{p}\right)}{P\left(\overline{L}^{(\emptyset)}_{m}=\overline{l}_m\mid \overline{A}_{p}=\overline{0},\overline{L}^{(\emptyset)}_{p}=\overline{l}_{p}\right)}\\
&=&\frac{P\left(Y_{k}^{(m)}=y,\overline{L}^{(\emptyset)}_{m}=\overline{l}_m\mid \overline{A}_{p-1}=\overline{0},\overline{L}^{(\emptyset)}_{p}=\overline{l}_{p}\right)}{P\left(\overline{L}^{(\emptyset)}_{m}=\overline{l}_m\mid \overline{A}_{p-1}=\overline{0},\overline{L}^{(\emptyset)}_{p}=\overline{l}_{p}\right)}\\
&=&\frac{P\left(Y_{k}^{(m)}=y,\overline{L}^{(\emptyset)}_{m}=\overline{l}_m,\overline{A}_{p-1}=\overline{0}\right)}{P\left(\overline{L}^{(\emptyset)}_{m}=\overline{l}_m,\overline{A}_{p-1}=\overline{0}\right)}\\
&=&P\left(Y_{k}^{(m)}=y|\overline{L}^{(\emptyset)}_{m}=\overline{l}_m,\overline{A}_{p-1}=\overline{0}\right),
\end{eqnarray*}
where for the first equality we used the induction hypothesis and for the fourth equality we used Assumption~\ref{intconf2}. This completes the induction step, and equation~(\ref{coarsealt}) and hence~(\ref{gammaalt}) follow.


\end{document}